\documentclass[12pt, a4paper]{article}

\usepackage{inputenc}
\usepackage{graphicx}
\usepackage{amssymb}
\usepackage{mathtools}
\usepackage{amsmath}
\usepackage{comment}
\usepackage{caption}
\usepackage{subcaption}
\usepackage{algorithm}
\usepackage{algorithmic}

\newcommand{\R}{\mathbb R}

\newtheorem{definition}{Definition}

\newtheorem{remark}{Remark}

\newtheorem{lemma}{Lemma}

\graphicspath{{./Figures/}}

\begin{document}
\begin{center}
{\large\bf{
The extension of linear inequality method for generalised rational Chebyshev approximation to approximation by general quasilinear functions.}  }
\end{center}

 \begin{center}{Vinesha Peiris and Nadezda Sukhorukova }
\end{center}


\begin{abstract}
In this paper we demonstrate that a well known linear inequality method developed for rational Chebyshev approximation  is equivalent to the application of the bisection method used in quasiconvex optimisation. Although this correspondence is not surprising, it naturally connects rational and generalised rational Chebyshev approximation problems with modern developments in the area of  quasiconvex functions and therefore offers more theoretical and computational tools for solving this problem. 
The second important contribution of this paper is the extension of the linear inequality method to a broader class of Chebyshev approximation problems, where the corresponding objective functions remain quasiconvex.  In this broader class of functions, the inequalities are no longer required to be linear: it is enough for each inequality to define a convex set and the computational challenge is in solving the corresponding convex feasibility problems. Therefore, we propose a more systematic and general approach for treating Chebyshev approximation problems. In particular, we are looking at the problems where the approximations are quasilinear functions with respect to their parameters, that are also the decision variables in the corresponding optimisation problems.
\end{abstract}


\section{Introduction}

In Chebyshev (uniform) approximation problems, the goal is to minimise the maximal deviation of the approximated function from the original function. In other words, the best uniform aproximation ensures that the deviation is minimal. The uniform polynomial approximation still draws an immense interest, since polynomials are simple and easy to handle. However, when nonsmooth or non-Lipschitz functions are subjected to approximate, polynomial approximation is not very efficient. 

Sometimes, this problem can be resolved by using piecewise polynomials or polynomial splines~\cite{NR,nurnberger,Schumaker1968,sukhorukovaoptimalityfixed}. Polynomial splines are more flexible  approximation technique than polynomials. Overall, polynomial spline approximation is efficient when the location of the knots (points of switching from one polynomial to another) is known. In the case of polynomials and piecewise polynomials with fixed knots, the optimisation problem is convex and therefore, existing optimisation techniques can be used to tackle this problem from the point of view of optimisation. When the knots are free (the location of the knots are unknown), the problem becomes nonconvex and the complexity of the corresponding optimisation problem increases~\cite{SukUgonCrouzeix2019,Meinardus1989,NR,nurnberger,Mulansky92,SukUgonTrans2017,su10}. In 1996, N\"urnberger, Daubechies, Borwein, and Totik~\cite{FreeKnotsOpenProblem96} identified this problem as a very hard and important open problem in approximation. The existing optimisation tools are not adapted to this problem due to its nonconvex and nonsmooth nature. One way to overcome this problem is to use rational functions (ratios of two polynomials).

Rational approximation attracted a significant research interest  in 50s-60s of the twentieth century~\cite{Achieser1965,Boehm1964,Meinardus1967rational,Ralston1965Reme,Rivlin1962} as a promising alternative to the free knot spline approximation. Moreover, there are theoretical studies~\cite{lorentz1996constructive,shevshevPopov1987}  demonstrated that approximation by rational functions and free-knot spline approximations are closely related to each other.

There are a number of efficient numerical methods for constructing rational approximations:  Remez method, originally developed for polynomial approximation, was extended to rational  approximation~\cite{fraser1962,ralston2001,rice1969}, the differential correction algorithm~\cite{DiffCorrection1972}, Loeb's algorithm~\cite{loeb1957rational}, the linear inequality method~\cite{loeb1960} and Osborne-Watson algorithm~\cite{osborne1969algorithm}. The most recent developments in this area \cite{Trefethen2018} are dedicated to ``nearly optimal'' solutions, whose construction is based on Chebyshev polynomials.  

The linear inequality method is interesting from several points of view. It is simple, easy to implement and converges to an optimal solution. Another important characteristic of this methods is its applicability to  generalised rational approximation problems~\cite{cheney1964generalized}, where approximations are constructed as the ratios of linear forms and not limited to ratios of polynomials. It was also demonstrated in~\cite{cheney1964generalized} that the corresponding optimisation problems in the case of generalised rational approximation are quasiconvex and strictly quasiconvex~\cite{barrodale1973} in the case of classical rational approximation. These results were used in the development of new theoretical characteristics of the problems rather than constructing computational methods. The first main purpose of this paper is to make a step towards covering this gap, in particular, we will demonstrate that the linear inequality method for generalised rational approximation is identical to the application of the bisection method developed for quasiconvex optimisation problems. The second main contribution of this paper is the extension of the linear inequality method to a broader class of Chebyshev approximation problems. This class includes all the approximation types, where the approximation is quasilinear with respect to its parameters. The corresponding inequalities may be nonlinear and therefore we change the name of the method to just ``inequality method'', but the set defined by each constraint is convex (that is, solving convex feasibility problems).

The rest of the paper is organised as follows. Section~\ref{sec:RationalApproximation} offers a detailed introduction on rational approximation and the linear inequality method. In section~\ref{sec:quasiconvexity}, we discuss the quasiconvexity of the problem and possible extension to a broader class of functions. Section~\ref{sec:bisection} provides of bisection method and its applicability to Chebyshev approximation problems. Section~\ref{sec:examples} discusses numerical experiments of the proposed method with the similarities and differences of the two methods while section~\ref{sec:conclusions} provides conclusion and future research directions.

\section{Rational and generalised  rational approximation}\label{sec:RationalApproximation}

\subsection{Problem formulation}\label{ssec:formulation}

Let $f(t)$ be a given real valued function defined on the closed interval $[c,d]$. Let $R_{n,m}$ be a rational function (rational function of type $(n,m)$) of degree $n$ in the numerator and degree $m$ in the denominator where $n$ and $m$ being nonnegative integers. 

$$R_{n,m}(t) =  P(t)/Q(t) = \frac{\sum_{i=0}^{n} p_i t^i}{\sum_{j=0}^{m} q_j t^j} .$$

The best uniform rational approximation problem over the closed interval $[c,d]$ is defined as follows:

\begin{equation} \label{eqn:best_rational_app}
\min_{r \in R_{n,m}} \max_{t \in [c,d]} \left | {f(t) - r(t)} \right |. 
\end{equation}
subject to,
\begin{equation*}
Q(t) > 0, t\in [c,d].
\end{equation*}

The objective is to determine the coefficients, $p_i$, $i = 0, \cdots ,n$ and $q_j$, $j = 0, \cdots ,m$ which minimise the maximal deviation. This problem~\eqref{eqn:best_rational_app} is also known as the minimax rational approximation problem. The rational function, $r(t)$ is defined to be a ratio of two polynomials. Most prominent methods (Remez method, differential correction method, linear inequality method, etc.), developed for uniform rational approximation were originally defined for the ratio of polynomials. Some of them (for example, linear inequality and differential correction) were extended to generalised  rational Chebyshev approximation (ratios of linear forms, introduced in~\cite{cheney1964generalized}).

In this paper we use the same definition of generalised rational Chebyshev approximation.
Specifically, we construct the approximations in the form of the ratios of linear combinations of functions in both numerator and the denominator.  The generalised rational approximation problem in Chebyshev norm can be formulated as follows~\cite{cheney1964generalized}:
\begin{equation} \label{eq:rational_app_obje_fun}
\min_{\substack{{\bf A}, {\bf B}}} \max_{t\in [c,d]} \left | f(t)-\frac{{\bf A}^T {\bf G}(t)}{{\bf B}^T {\bf H}(t)} \right |. 
\end{equation}\\
subject to,\\
\begin{equation}
{{\bf B}^T {\bf H}(t)} > 0, t\in [c,d]. \label{eq:rational_app_const}
\end{equation}\\
where,

\begin{itemize}
	\item $f(t)\in {C_{[c,d]}^0}$ is the function to approximate,
	\item ${\bf A} = (a_0,a_1,\cdots,a_n)^T \in \mathbb{R}^{n+1}$ and ${\bf B} = (b_0,b_1,\cdots,b_m)^T \in \mathbb{R}^{m+1}$ are the decision variables,
	\item ${\bf G}(t) = (g_0(t),g_1(t),\cdots,g_n(t))^T$ and ${\bf H}(t) = (h_0(t),h_1(t),\cdots,h_m(t))^T$, $g_j, j=1,\cdots,n$ and  $h_i, i=1,\cdots,m$ are known functions which are also known as the basis functions.
\end{itemize}

The constraint set (strict linear inequality) is an open convex set. The numerator and  the denominator are linear combinations of basis functions. When all the basis functions are monomials, the approximations are rational functions from $R_{n,m}$, and the problem is reduced to the best uniform rational approximation~\eqref{eqn:best_rational_app}. Our approximations are not restricted to rational functions with polynomials. All of the results are valid for any type of basis functions when~\eqref{eq:rational_app_const} is satisfied.
\subsection{Linear inequality method}\label{ssec:inequality} 

The linear inequality method~\cite{loeb1960} for solving rational and generalised rational problems includes the following steps.
\begin{enumerate}
\item[Step 0] Identify upper and lower bounds ($M_u$ and $M_l$) for the maximal deviation. In particular, zero is always a  lower bound, while $$\max_{t\in[a,b]}f(t)-\min_{t\in[a,b]}f(t)$$ can be used as an upper bound. 
\item[Step 1]  Set $z=\frac{1}{2}(M_u+M_l)$. Check if the following system of inequalities has a feasible solution.
\begin{equation}
f(t)-\frac{{\bf A}^T {\bf G}(t)}{{\bf B}^T {\bf H}(t)}\leq z, \label{eq:rational_app_feas1}
\end{equation}
\begin{equation}
\frac{{\bf A}^T {\bf G}(t)}{{\bf B}^T {\bf H}(t)}-f(t)\leq z, \label{eq:rational_app_feas2}
\end{equation}
\begin{equation}
{{\bf B}^T {\bf H}(t)} >0, t\in [c,d]. \label{eq:rational_app_feas}
\end{equation}
\item[Step 2] If (\ref{eq:rational_app_feas1})- (\ref{eq:rational_app_feas}) has a feasible solution, assign $M_u=z$, otherwise $M_l=z$. If  $M_u-M_l\geq \varepsilon$, go to Step~1.
\end{enumerate}
Therefore, the algorithm stops when  $M_u-M_l<\varepsilon$ (tolerance).

It was demonstrated in~\cite{cheney1964generalized} that in the case when $t$ takes discrete values (discretised system), Step~1 can be verified by solving a linear programming problem. For computations purposes, constraint~(\ref{eq:rational_app_feas}) is replaced by 
\begin{equation}
{{\bf B}^T {\bf H}(t)} \geq \delta, t\in [c,d], . \label{eq:rational_app_feas0}
\end{equation}
where $\delta$ is a small positive number.
%
%
%
%

\section{Quasiconvexity of the problem}\label{sec:quasiconvexity}
\subsection{Quasiconvex functions}\label{ssec:quasiconvexity}

%

The notion of quasiconvexity was first introduced by mathematicians working with the area of financial mathematics. 
There are two commonly used equivalent definitions of quasiconvex functions.
\begin{definition} \label{def:quasiconvex1}
	A function $ f : D \rightarrow \mathbb{R} $ defined on a convex subset $D $ of a real vector space is called  quasiconvex if and only if  for any pair $x$ and $y$ from $D$  and $\lambda\in [0,1]$ one has 
	$$ 
	f (\lambda x + ( 1-\lambda ) y )
	\leq
	\max\{ f ( x ) , f ( y )\} .$$ 
\end{definition}
\begin{definition}  \label{def:quasiconvex2}
	Function $ f ( t )$ is quasiconvex if and only if its  sublevel set 
	$$ S_{\alpha}=\{x |f (x ) \leq \alpha \} $$
	is  convex for any $\alpha \in \mathbb{R}$.  The set $S_{\alpha}$ is also called $\alpha$-sublevel set.
\end{definition}
Definition~\ref{def:quasiconvex1} is due to to Bruno de Finetti's work~\cite{deFinetti1949}. In this paper de Finetti studied the behaviour of functions, whose sublevel set is convex. The term quasiconvexity was introduced much later, as a result to generalised the notion of convexity to a broader class of functions.

\begin{definition} \label{def:quasiconcave}
	Function $f$ is quasiconcave if and only if $-f$ is quasiconvex.
\end{definition}
This means that every superlevel set $ \bar{S}_{\alpha}=\{x |f (x ) \geq \alpha \} $ is convex.
\begin{definition} \label{def:quasiaffine}
	Functions that are quasiconvex and quasiconcave at the same time are called quasiaffine (sometimes quasilinear).
\end{definition}

There are many studies devoted to quasiconvex optimisation~\cite{JPCrouzeix1980quasi,DaniilidisHadjisavvasMartinezLegas2002,dutta2005abstract,RubinovSimsek,Rubinov00} just to name a few.

It was demonstrated in~\cite{cheney1964generalized} that the sublevel sets of the objective functions~\eqref{eq:rational_app_obje_fun} are convex (\cite{cheney1964generalized}, Lemma~2). This result was considered by the authors as an auxiliary result, but this also proves that the objective functions~\eqref{eq:rational_app_obje_fun} is quasiconvex. The same result can be obtained by applying the following two properties of quasiconvex functions. 
\begin{enumerate}
	\item The supremum of a family of quasiconvex functions is quasiconvex.
	\item The ratio of two affine functions is quasiaffine.
\end{enumerate}
\begin{remark}
In the case of classical rational approximation (that is, all the basis functions are monomials), it was proved that the  objective function is strictly quasiconvex~\cite{barrodale1973}.  
 \end{remark}
Therefore, the problem of rational and generalised rational approximation can be treated using a number of computational methods developed for quasiconvex optimisation~\cite{MLegazquasiconvexduality,DaCruzAlgorithmsQuasiconvex}. One such method (called bisection method for quasiconvex functions) will be described in section~\ref{sec:bisection}. In the same section we demonstrate that the bisection method for quasiconvex functions is equivalent to the linear inequality method for classical rational approximation and can be extended to generalised rational approximation (the constraints remain linear) or even to a broader class approximations, where the corresponding optimisation problems remain quasiconvex. This class of approximations is described in the next section.

\subsection{Extension to broader classes of functions}\label{ssec:extension}

The goal of this section is to identify the class of Chebyshev approximations whose optimisation problems are quasiconvex and therefore can be solved by applying quasiconvex optimisation techniques. The following lemma holds.

\begin{lemma}\label{lemma:class_of_approximations}
If the approximation is quasilinear with respect to its parameters, that are also the decision variable for the corresponding optimisation problems, then the Chebyshev approximation problem is quasiconvex.  
\end{lemma}
{\bf Proof:} We need to prove that if the approximations are quasilinear with respect to their parameters then the corresponding optimisation problem is quasiconvex. The optimisation problem is as follows:
\begin{equation}\label{eq:approx_quasilinear}
\min_A\max_{t\in[c,d]}|F(A,t)-f(t)|,
\end{equation}
where $f(t)$ is the function to approximate, $F(A,t)$ is the approximation and $A$ are the parameters of approximation that can be optimised. We need to prove that if $F(A,t)$ is quasilinear with respect to $A$ then $$\Phi(A)=\max_{t\in[c,d]}|F(A,t)-f(t)|$$ is quasiconvex. Indeed, 
$$\Phi(A)=\max_{t\in[c,d]}|F(A,t)-f(t)|=\max_{t\in[c,d]}\{F(A,t)-f(t),f(t)-F(A,t)\}.$$
Since $F(A,t)$ is quasilinear for $A$, both functions $F(A,t)-f(t)$ and $f(t)-F(A,t)$ are quasiconvex.  Then the maximum of quasiconvex functions is also quasiconvex.

\hskip 300pt
$\square$

\section{Bisection method for quasiconvex functions}\label{sec:bisection}

In this section we describe a simple approach for solving quasiconvex problems. This method is also known as a bisection method for quasiconvex functions (see~\cite{SL}, section~4.2.5 for more information). The approach is based on the representation of the sublevel set of quasiconvex functions. These sets are convex (due to the definition of quasiconvex functions). 

Let $\varphi_z(x): \R^k\rightarrow\R$ and $z\in\R$ be a family of convex functions, such that $\varphi_z(x)\leq 0$ if and only if $x$~belongs to $z$-sublevel set of the quasiconvex function $\Phi$, that is equivalent to $\Phi(x)\leq z$.

Consider the following quasiconvex optimiation problem:
\begin{equation} \label{eq:quasi_obje_fun}
\min_x \Phi(x), 
\end{equation}
subject to
\begin{equation}
\Phi_i(x)\leq 0, i=1,\dots,l, \label{eq:quasi_const}\\
\end{equation}
\begin{equation}
Ax=b\label{eq:quasi_const1},
\end{equation}
where $\Phi(x)$ is a quasiconvex function, $\Phi_i,$ $i=1,\dots,l$ are convex functions and (\ref{eq:quasi_const1}) is a set of linear equations. Then the feasibility problem is
\begin{equation*}
{\rm Find}~ x,~{\rm subject~to} 
\end{equation*}
\begin{equation}
\varphi_z(x)\leq 0, \label{eq:feas_quasi_const0}\\
\end{equation}
\begin{equation}
\Phi_i(x)\leq 0, i=1,\dots,l, \label{eq:feas_quasi_const}\\
\end{equation}
\begin{equation}
Ax=b.\label{eq:feas_quasi_const1}
\end{equation}

Since all the sublevel sets are convex, the goal of the feasibility problem~(\ref{eq:feas_quasi_const0})-(\ref{eq:feas_quasi_const1}) is to find a point that belongs to a convex set (intersection of convex sets). This class of problems is also called convex feasibility problems. These problems may be very challenging, but there  are a number of efficient techniques for solving convex feasibility problems~(\cite{BauschkeLewis, ACNN, Shi-yaXu, NimitPetrot,  YangYang, Zaslavski, Zhao}~just to name a few). There are still several open problems. 
However, in this paper we do not discuss these problems and concentrate of the approximation side.  

In the case of classical rational and generalised rational approximation, the corresponding convex feasibility problems (discrete case) can be reduced to solving linear programming problems~\cite{cheney1964generalized,loeb1960,loeb1957rational,osborne1969algorithm,AMCPeirisSukhSharonUgon}, while for the continuous case the corresponding convex feasibility problems  are  linear semi-infinite. More details on linear semi-infinite problems can be found in~\cite{goberna1998comprehensive}.

In the case of discrete approximation by rational and generalised rational functions, for a given level $\alpha=z$, the corresponding linear programming problem is as follows:
\begin{equation}\label{eq:obja}
\min u
\end{equation}
subject to
\begin{align}
(f(t_i)-z){{\bf{B}}^T{\bf{H}}(t_i)}-{{\bf{A}}^T{\bf{G}}(t_i)}\leq u,
\quad i=1,\dots,N  \label{eq:con1a} \\
{{\bf{A}}^T{\bf{G}}(t_i)}-(f(t_i)+z){{\bf{B}}^T{\bf{H}}(t_i)}\leq u, \quad i=1,\dots,N  \label{eq:con2a} \\
{\bf{B}}^T{\bf{H}}(t_i)\geq\delta,
\quad i=1,\dots,N,  \label{eq:con3a} 
\end{align}
where $\delta$ is a small positive constant, sublevel set parameter~$z$ is also a constant and therefore (\ref{eq:obja})-(\ref{eq:con3a}) is a linear programming problem.

In the case of discrete approximation by a general quasiconvex function with respect to the parameters, for a given level $\alpha=z$, the corresponding linear programming problem is as follows:
\begin{equation}\label{eq:objaquasi}
{\rm Find}~A
\end{equation}
subject to
\begin{align}
f(t_i)-F(A,t_i)\leq z,
\quad i=1,\dots,N  \label{eq:con1aquasi} \\
F(A,t_i)-f(t_i)\leq z, \quad i=1,\dots,N.  \label{eq:con2aquasi} \\
\end{align}
Since $F(A,t)$ is quasilinear, the feasible set described by inequalities~(\ref{eq:con1aquasi})-(\ref{eq:con2aquasi}) is convex.

We use a bisection method for quasiconvex optimisation. This method is simple and reliable and can be applied to any quasiconvex function (see~\cite{SL}, section~4.2.5 for more information). Our bisection method relies on solving the  corresponding convex feasibility problem,  which can be done by solving~(\ref{eq:obja})-(\ref{eq:con3a}) in the case of rational and generalised rational approximation. The bisection method is given in Algorithm~\ref{alg:bisection}.
\begin{algorithm}
	\caption{Bisection algorithm for quasiconvex optimisation} \label{alg:bisection}
	\begin{algorithmic}
		\ENSURE Maximal deviation $z$ (within $\varepsilon$ precision). 
Start with given precision $\varepsilon>0$.
		\STATE set $l \leftarrow 0$
		\STATE set $u$
		\STATE $z \leftarrow (u+l)/2$
		\WHILE{$u-l \leq \varepsilon$}
		\STATE {Check if problem~\eqref{eq:objaquasi}-\eqref{eq:con2aquasi} has a feasible solution with $z$.}
		 \IF{feasible solution exists}
		\STATE $u \leftarrow z$
		\ELSE
		\STATE $l \leftarrow z$
		\ENDIF
		\STATE update $z \leftarrow (u+l)/2$
		\ENDWHILE 
	\end{algorithmic}
\end{algorithm}

Therefore, the linear inequality method developed for rational and generalised rational approximation is equivalent to applying the bisection method for quasiconvex optimisation. In the case, the convex feasibility problem can be formulated as a linear programming problem. In more general cases, where the approximations are arbitrary quasilinear functions with respect to their parameters (Lemma~\ref{lemma:class_of_approximations}) more general convex feasibility problems have to be solved. In this case, the set of inequalities are not restricted to linear: the only requirement for the left hand side is to be quasilinear.

\section{Numerical examples}\label{sec:examples}
\subsection{Free knot linear approximations}\label{ssec:free_knot}

We start with a numerical example from the area of free knot spline approximation. A continuous function $f(t)$ is to be approximated by a piecewise linear function with one internal knot (point of switching from one linear piece to another). In general, in the case free knot polynomial spline approximation the corresponding optimisation problems are very difficult. When the knots are fixed and only the coefficients of polynomial pieces are subject to approximation, the problems are convex. Therefore, there are a number of methods that include block-coordinate approximation~\cite{nurnberger}, where the polynomial coefficients are subject to optimisation (knots are fixed) and then the obtained coefficients are fixed and the knots are subject to optimisation. 

Suppose that a continuous function $f(t)$ is approximated on an interval $[c,d]$ by a piecewise linear function 
$$S(a_0,a_1,a_2,\theta,t)=a_0+a_1t+a_2\max\{0, t-\theta\}.$$
Assume that $a_0, a_1$ and $a_2$ (coefficients) are fixed and  the knot $\theta$ is subject to optimisation. Then the corresponding optimisation problem is as follows:
\begin{equation}\label{eq:exfreeknotlinprob}
\min_{\theta\in[c,d]}\max_{t\in[c,d]}|f(t)-(a_0+a_1t+a_2\max\{0,t-\theta\})|.
\end{equation} 
Assume that $a_2\neq 0$, otherwise it is equivalent to $\theta=d$ (that is, only one linear piece).
Therefore, there is only one variable $\theta$.  

The function $\max\{0,t-\theta\}$ is quasilinear with respect to $\theta$ and therefore $\Phi(\theta)=\max_{t\in[c,d]}|f(t)-(a_0+a_1t+a_2\max\{0,t-\theta\})|$ is quasiconvex in $[c,d]$. Since there is only one variable, there are several methods to minimise $\Phi(\theta)$. We use this example for illustrative purposes to interpret bisection.

Fix $\alpha=\Delta\geq 0$ and then one needs to check if there exists a value $\theta=\theta_0$, such that
\begin{equation}\label{eq:estimation}
f(t)-a_0-a_1t-\Delta\leq a_2\max_{t\in[c,d]}\{0,t-\theta\}\leq f(t)-a_0-a_1t+\Delta,
\end{equation}
then
\begin{equation}\label{eq:estimation1}
g(t)-z\leq\max_{t\in[c,d]}\{0,t-\theta\}\leq g(t)+z,
\end{equation}
where $g(t)= \frac{f(t)-a_0-a_1t}{a_2}$, $z=\Delta/|a_2|$ and therefore $z\geq 0$. Hence, the feasibility problem is equivalent to fitting $\max\{0,t-\theta\}$ in the area between $g(t)+z$ and $g(t)-z$. 

\begin{figure}
  \center
  \includegraphics[width=0.9\linewidth]{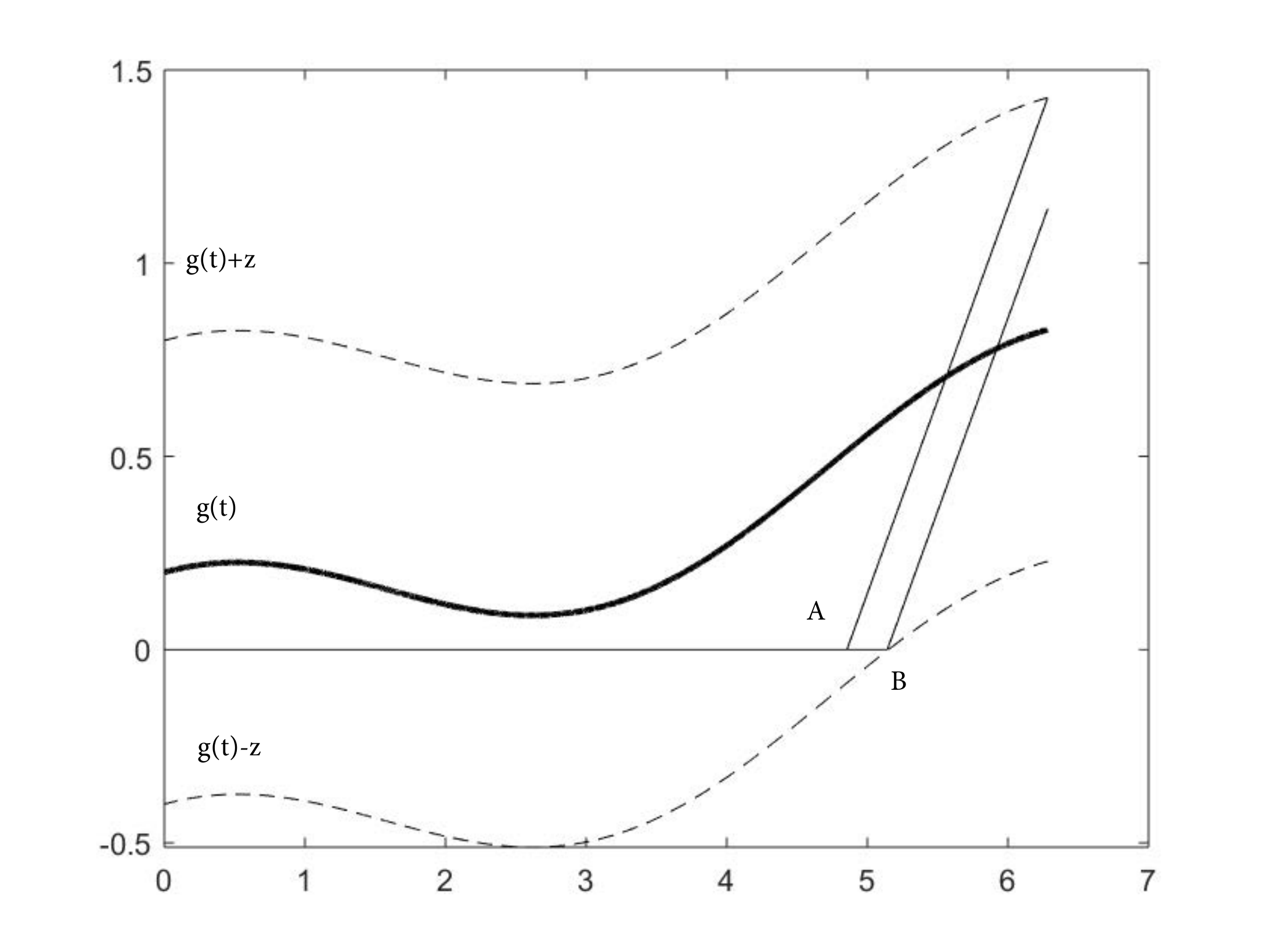}
  \caption{More than one possibility for $\theta$, any value between $A$ and $B$ gives a feasible solution.}
  \label{fig:severalsolutions}
\end{figure}

\begin{figure}
  \center
  \includegraphics[width=0.9\linewidth]{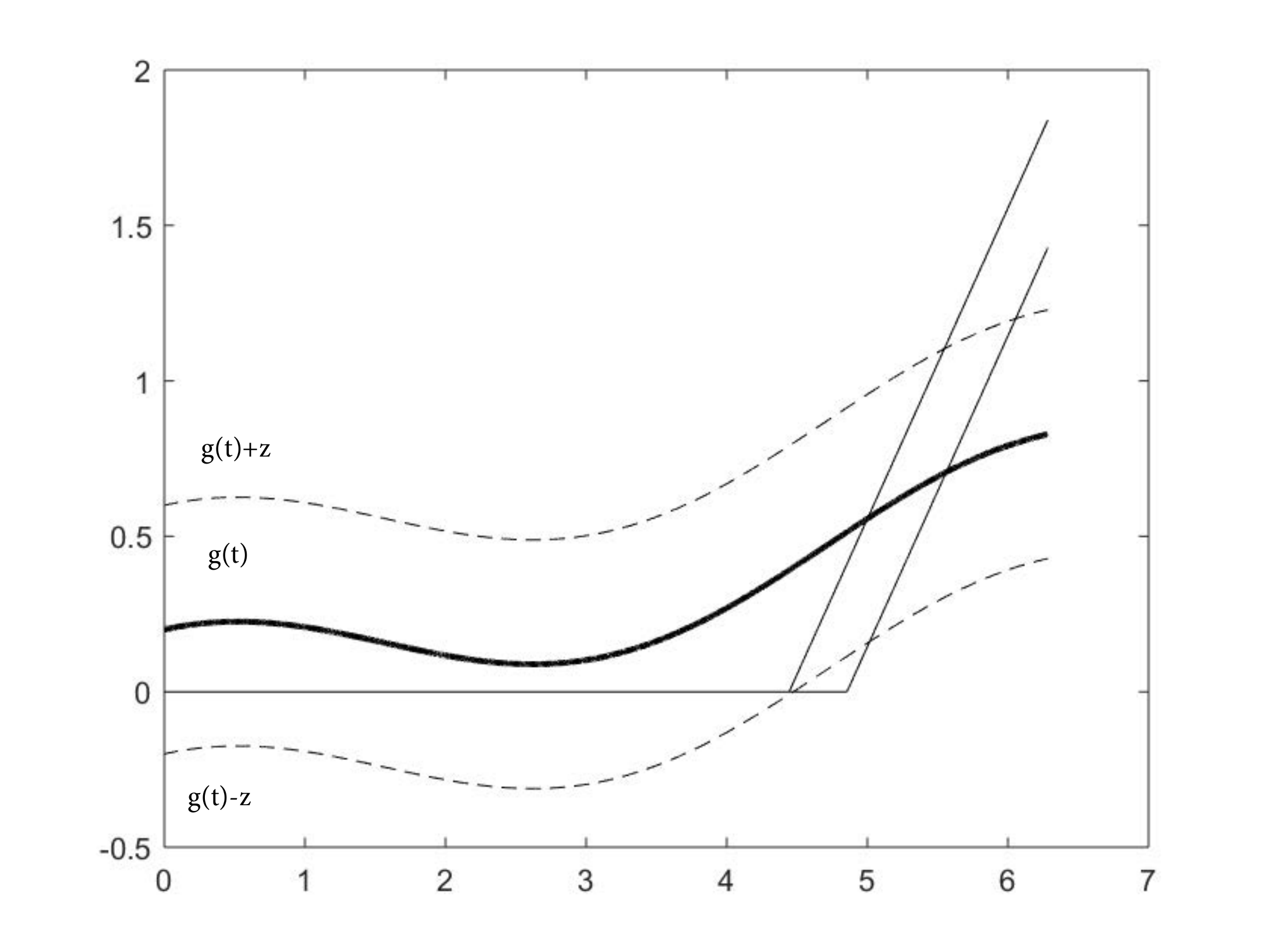}
  \caption{There is no $\theta$, not feasible solution.}
  \label{fig:nosolution}
\end{figure}

Figure~\ref{fig:severalsolutions} depicts the situation there is more than one value for $\theta$, such that there is a function of the form $\max\{0,t-\theta\}$ that fits in the area between $g(t)+z$ and $g(t)-z$. Indeed, if $A=(t_1,0)$ and $B=(t_2,0)$, then $t_1\leq\theta\leq t_2$. Figure~\ref{fig:nosolution} depicts the situation where the corresponding feasibility problem has no solution. In this case, there is no function of the form $\max\{0,t-\theta\}$ that fits in the area between $g(t)+z$ and $g(t)-z$. This means that there is no~$\theta$ that satisfies the corresponding feasibility problem.

\subsection{Ratio of piecewise polynomial approximations}\label{ssec:piecewise_ratio}

In this section we present two sets of numerical experiments where the function $f(t)$ is approximated by a ratio of two piecewise polynomials. In both cases, the ratio consists of linear combinations of basis functions in the numerator and in the denominator. The piecewise polynomials are constructed using the truncated power function:
\begin{equation*}
	S({\bf A},t,\theta) = a_0 + \sum_{j=1}^{m}{a_j t^j} +\sum_{k=1}^{m} {a_{k+m} (\max\{0,t-\theta\})^k}
\end{equation*}
where $m$ is the degree of the piecewise polynomial, $a_i \in {\bf A}, i=0,\dots,2m+1$ are the parameters and assume that the knots are fixed.

We consider approximating a nonlinear, non-Lipschitz function
$$f(t) = \sqrt{|t-0.25|}, \quad t \in [-1,1].$$

The approximations are calculated via Algorithm~\ref{alg:bisection} and the precision ($\varepsilon$) is set to be $10^{-5}$. The code for these numerical experiments are implemented using Matlab.

\subsubsection{First set of experiments}\label{sssec:first_set}

Suppose that the function $f(t)$ is approximated by a ratio of two piecewise polynomials with fixed knots where the degrees of the numerator and the denominator are $2,2$ respectively.

The corresponding optimisation problem is as follows:
\begin{equation*}
\min_{A,B} \max_{t \in [-1,1]} \left| f(t) - \frac{a_0 + a_1t + a_2t^2 + a_3\max\{0,t-\theta_1\} + a_4(\max\{0,t-\theta_1\})^2}{1 + b_1t + b_2t^2 + b_3\max\{0,t-\theta_2\} + b_4(\max\{0,t-\theta_2\})^2} \right|
\end{equation*}
where, $a_i \in {\bf A}, i=0,\dots,4$ and $b_j \in {\bf B}, j=1,\dots,4.$

The set of parameters (${\bf A,B}$), are the coefficients of piecewise polynomials and there are $9$ parameters in total (computed as $5+4$) to be determined. Both numerator and denominator consist of linear functions with respect to the parameters. Therefore, this ratio forms a quasilinear function which leads to the bisections method defined for quasiconvex problems. Note that the constant term of the denominator is fixed at 1 to avoid the ambiguity of two solutions representing the same fraction.

One can fix the knots at anywhere between $[-1,1]$. Let us assume that both knots in the numerator and the denominator are located at the same point. That is, $\theta_1 = \theta_2 (=\theta)$ and the approximation comprise of two intervals $[-1,\theta], [\theta,1]$. We consider four different cases where~$\theta$ is located at 0.25, 0.5, 0 and -0.5. 

\begin{figure}
	\begin{center}
		\begin{subfigure}{.49\textwidth}
			\includegraphics[width=\textwidth]{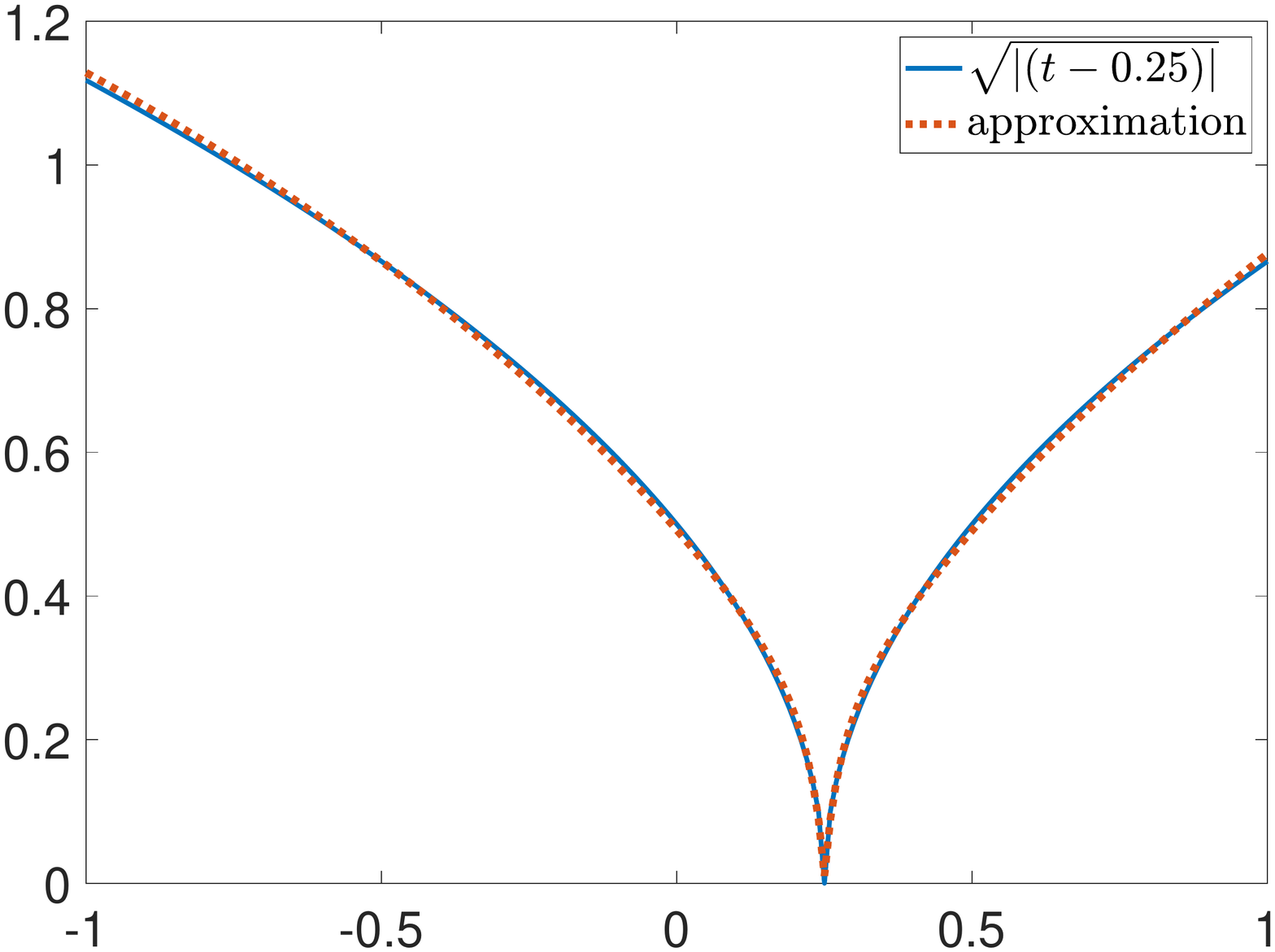}  
			\caption{The function and its approximation}
			\label{subfig:fun_app_theta_0.25}
		\end{subfigure}
		\begin{subfigure}{.49\textwidth}
			\includegraphics[width=\textwidth]{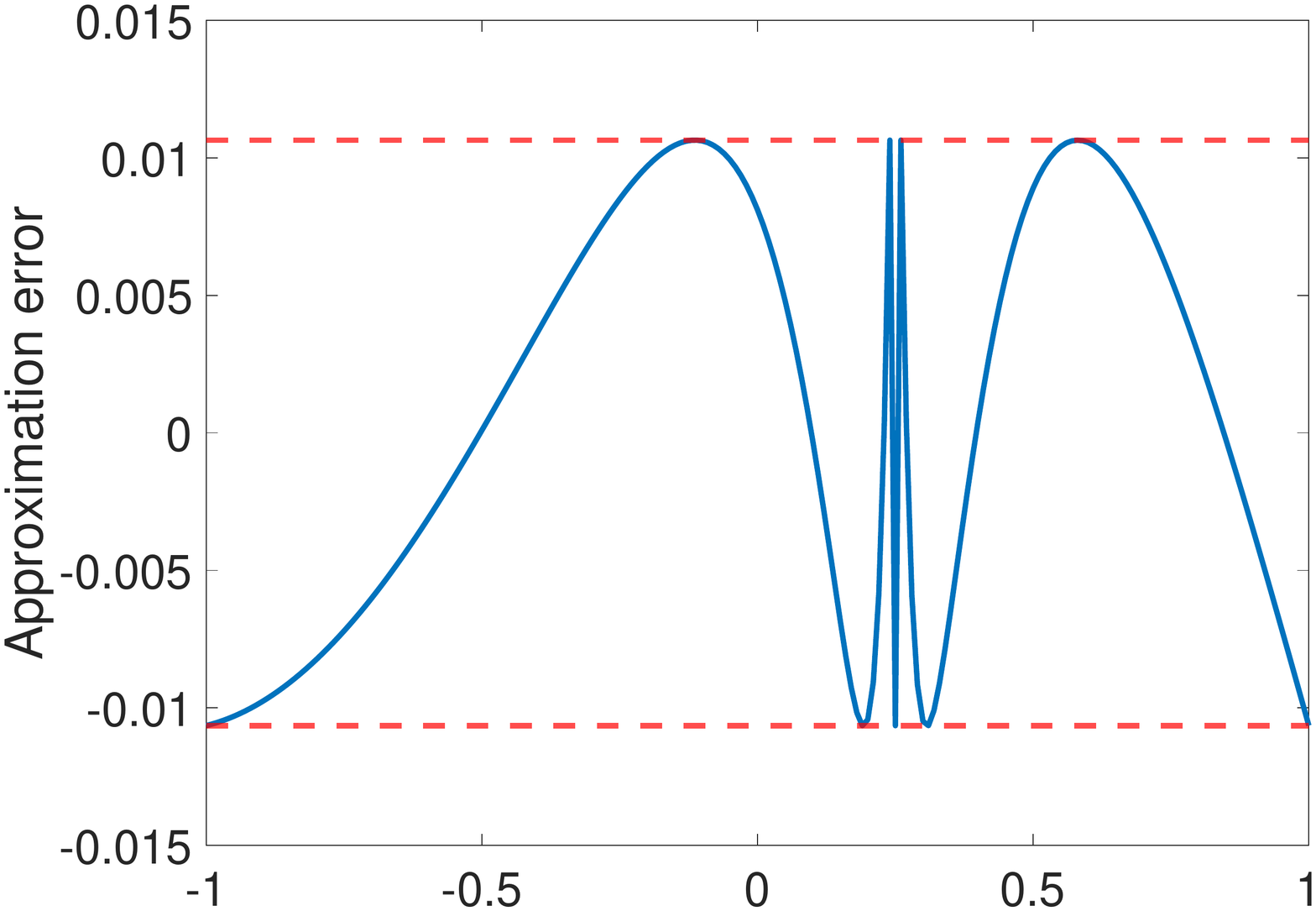}  
			\caption{Error curves}
			\label{subfig:error_curve_theta_0.25}
		\end{subfigure}
		\caption{The knot ($\theta_1=\theta_2$) is located at 0.25}
		\label{fig:same_theta_0.25}
	\end{center}
\end{figure}

Figure~\ref{fig:same_theta_0.25} shows that when $\theta = 0.25$, i.e. the knot is located at the same point where the original function tends to have an abrupt change,  the alternation sequence contains 9 (5+4) maximum deviation points, and hence the solution is optimal. The maximal deviation is around 0.01. In~\cite{AMCPeirisSukhSharonUgon} the same function $f(t)$ was approximated by a rational function where the basis functions are just monomials and the maximum deviation was around 0.05. Therefore, if the exact location of the knot is correctly identified, it is more likely to get an enhanced result by using appropriate basis functions which forms a quasilinear function.

\begin{figure}
	\begin{center}
		\begin{subfigure}{.49\textwidth}
			\includegraphics[width=\textwidth]{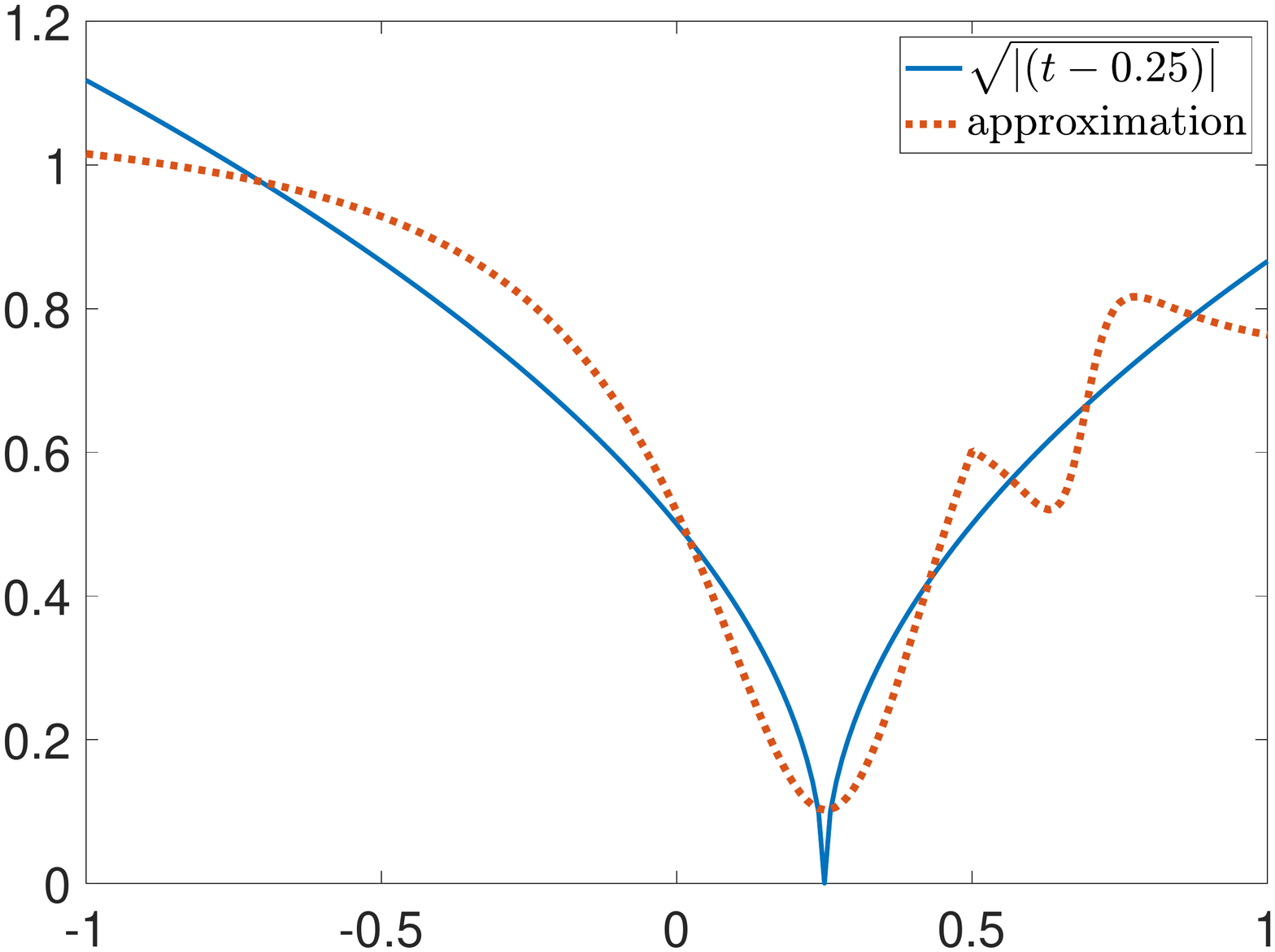}  
			\caption{The function and its approximation}
			\label{subfig:fun_app_theta_0.5}
		\end{subfigure}
		\begin{subfigure}{.49\textwidth}
			\includegraphics[width=\textwidth]{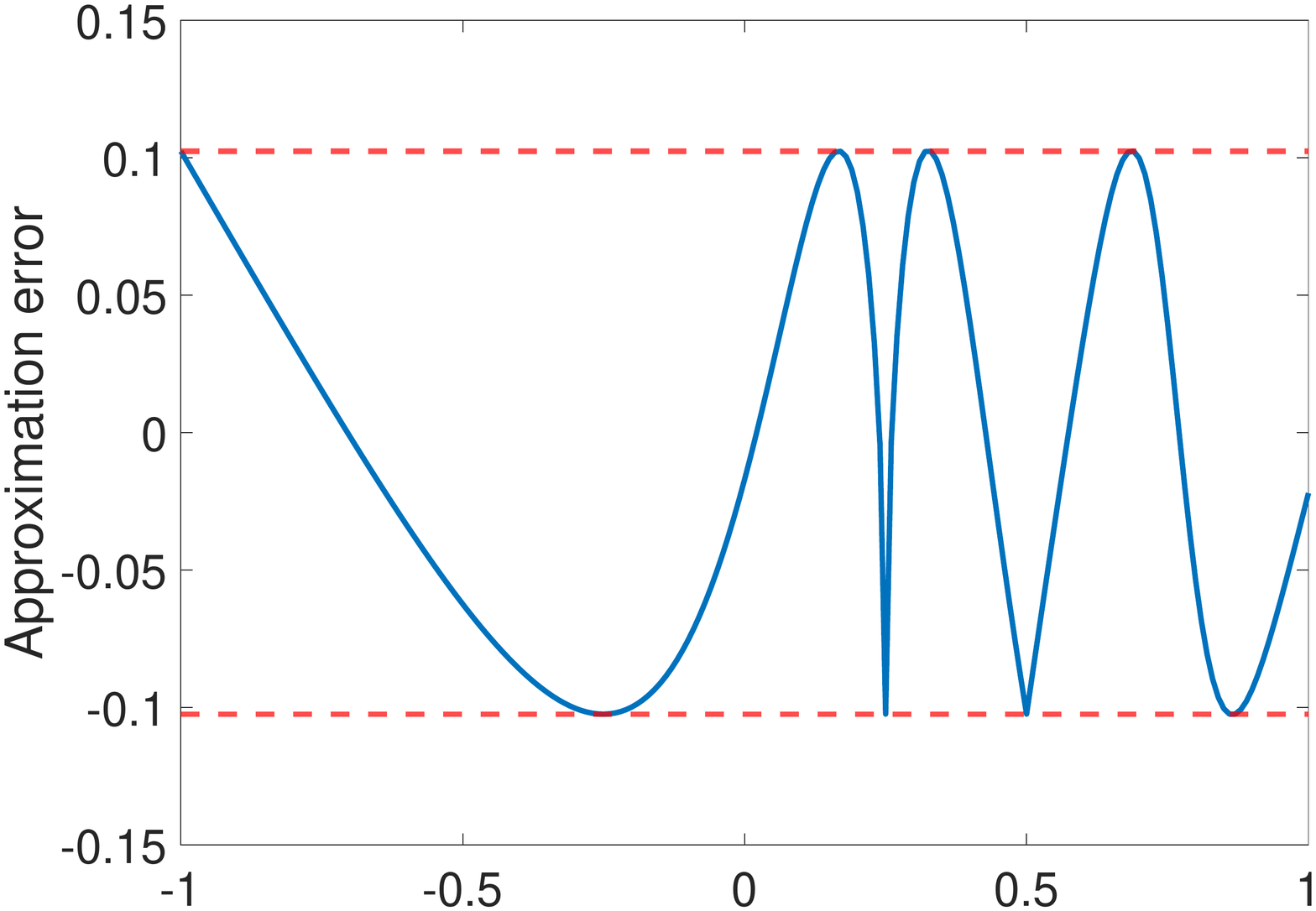}  
			\caption{Error curves}
			\label{subfig:error_curve_theta_0.5}
		\end{subfigure}
		\caption{The knot ($\theta_1=\theta_2$) is located at 0.5}
		\label{fig:same_theta_0.5}
	\end{center}
\end{figure}

Figure~\ref{fig:same_theta_0.5} illustrates 8 alternating error peaks that come with the same magnitude when the $\theta$ is located at 0.5. However, by decreasing the precision parameter of the algorithm~\ref{alg:bisection}, we may possibly end up with 9 alternating error peaks. The maximum deviation is around 0.1 which is 10 times higher than the previous error. 

\begin{figure}
	\begin{center}
		\begin{subfigure}{.49\textwidth}
			\includegraphics[width=\textwidth]{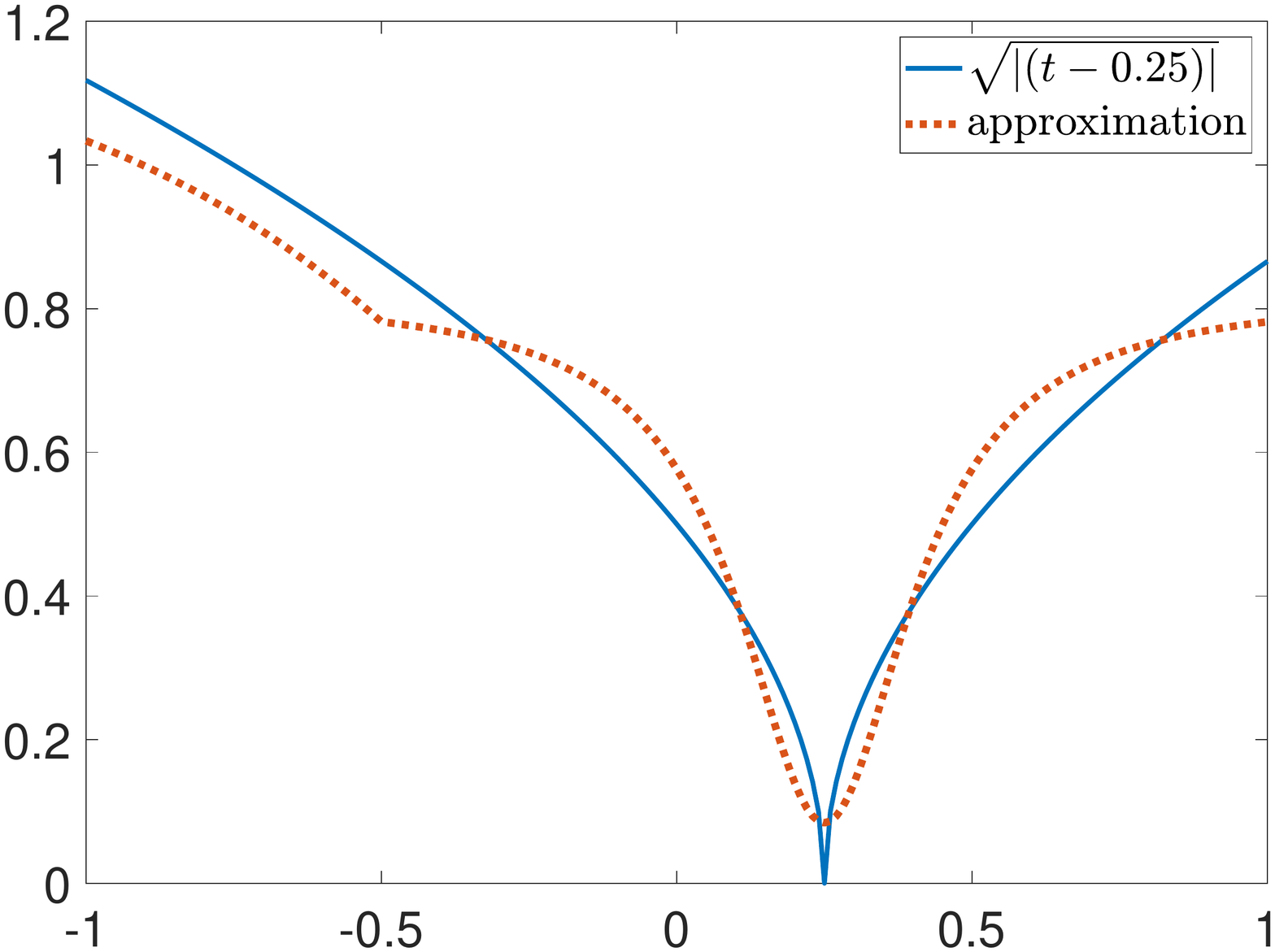}  
			\caption{The function and its approximation}
			\label{subfig:fun_app_theta_-0.5}
		\end{subfigure}
		\begin{subfigure}{.49\textwidth}
			\includegraphics[width=\textwidth]{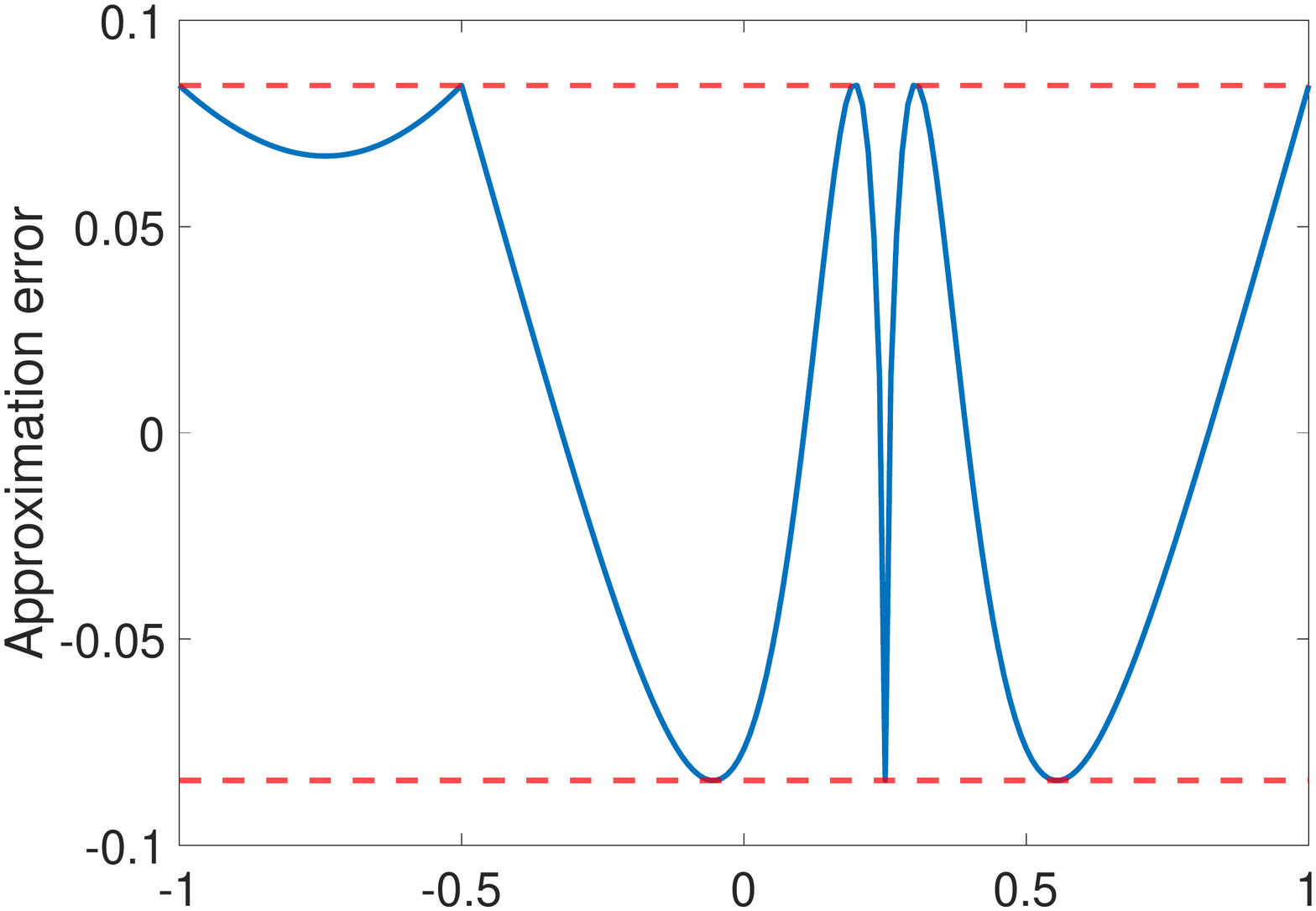}  
			\caption{Error curves}
			\label{subfig:error_curve_theta_-0.5}
		\end{subfigure}
		\caption{The knot ($\theta_1=\theta_2$) is located at -0.5}
		\label{fig:same_theta_-0.5}
	\end{center}
\end{figure}

\begin{figure}
	\begin{center}
		\begin{subfigure}{.49\textwidth}
			\includegraphics[width=\textwidth]{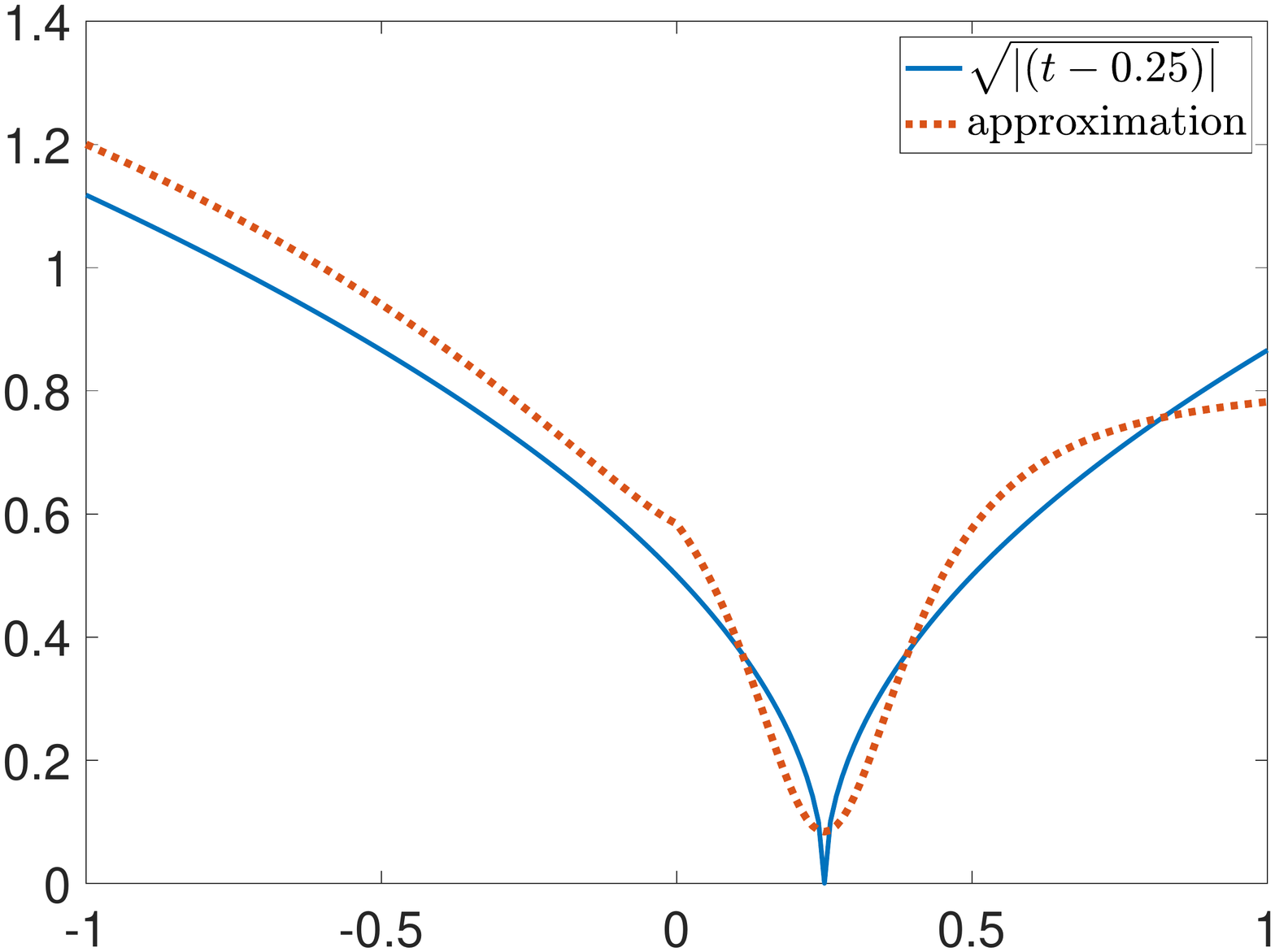}  
			\caption{The function and its approximation}
			\label{subfig:fun_app_theta_0}
		\end{subfigure}
		\begin{subfigure}{.49\textwidth}
			\includegraphics[width=\textwidth]{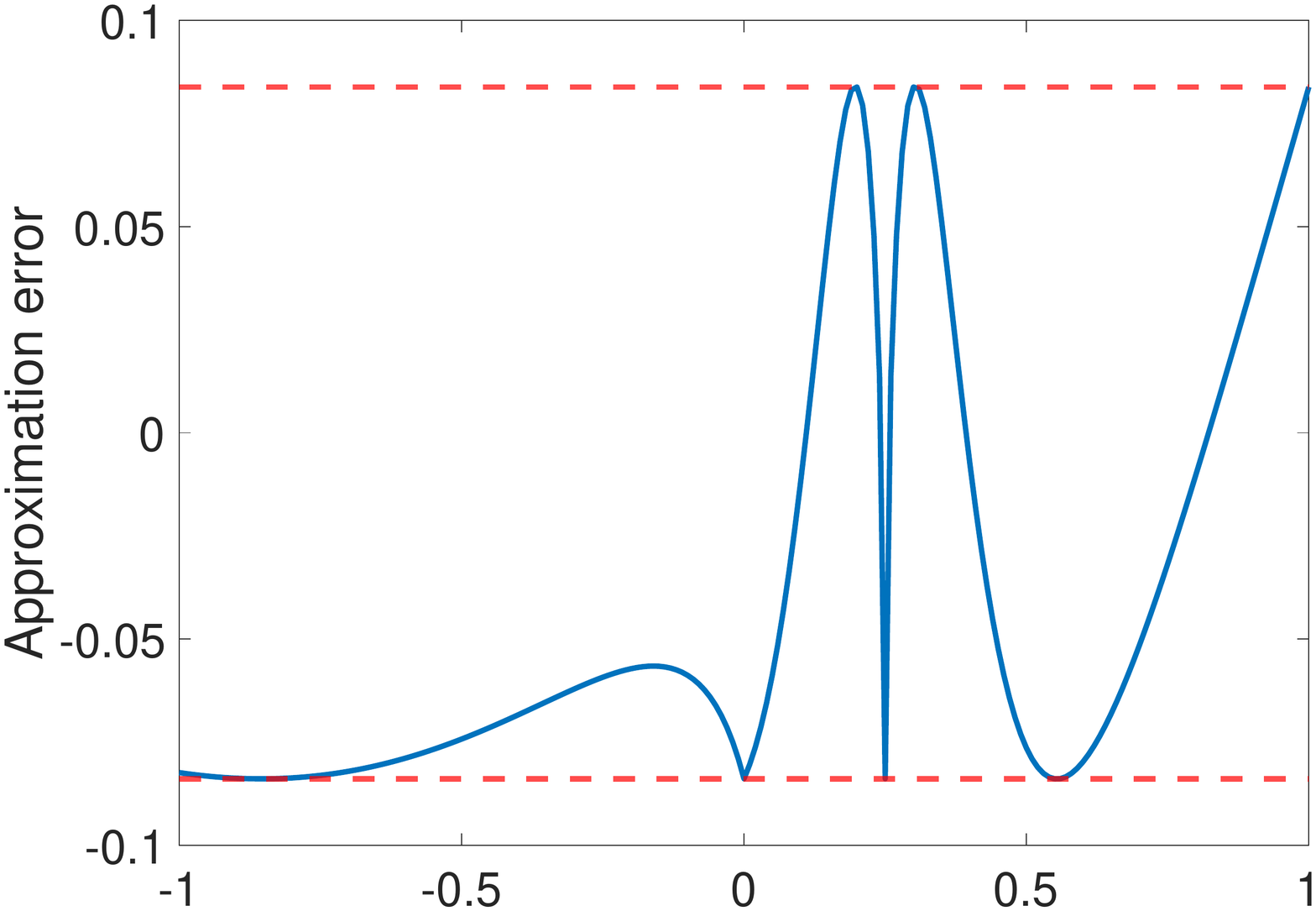}  
			\caption{Error curves}
			\label{subfig:error_curve_theta_0}
		\end{subfigure}
		\caption{The knot ($\theta_1=\theta_2$) is located at 0}
		\label{fig:same_theta_0}
	\end{center}
\end{figure}

Figure~\ref{fig:same_theta_-0.5} and Figure~\ref{fig:same_theta_0} depict the corresponding approximations and the error curves when $\theta = -0.5$ and $\theta = 0$ respectively. There are less than 9 maximal deviation and the magnitude of the error is around 0.1 for both cases.

Let us now assume that the location of the knot in the numerator is different to the location of the knot in the denominator. That is, the approximation consist of 3 intervals $[-1,\theta_1], [\theta_1,\theta_2]$ and $[\theta_2,1]$ or $\theta_1$ and $\theta_2$ are in different order. We consider different combinations for the locations of the knots as stated in table~\ref{tab:1h}.

\begin{table}
	\centering
	\caption{Combinations for the locations of the knots}
	\label{tab:1h}
	\begin{tabular}{ |r|r|c|c| } 
		\hline
		$\theta_1$ & $\theta_2$ &Number of deviation points & Figure\\
		\hline
		0.25            & 0.5             & 7    & \ref{subfig:error_curve_theta_0.25,0.5}\\
		0.25            & 0               & 7    & \ref{subfig:error_curve_theta_0.25,0}\\
		0.5             & 0.25            & 5    & \ref{subfig:error_curve_theta_0.5,0.25}\\
		0.5             & 0               & 6    & \ref{subfig:error_curve_theta_0.5,0}\\
		0               & 0.25            & 7    & \ref{subfig:error_curve_theta_0,0.25}\\
		0               & 0.5             & 7    & \ref{subfig:error_curve_theta_0,0.5}\\
		$-\frac{1}{3}$  & $\frac{1}{3}$   & 7    & \ref{subfig:error_curve_theta_(-1/3,1/3)}\\
		$\frac{1}{3}$   & $-\frac{1}{3}$  & 9    & \ref{subfig:error_curve_theta_(1/3,-1/3)}\\
		\hline
	\end{tabular}
\end{table}

In both Figure~\ref{fig:diff_theta_0.25,0.5} and Figure~\ref{fig:diff_theta_0.25,0}, the approximations appear to be very similar except for the location of the knot in the denominator. The knot in the numerator is identified and located correctly. Both error curves have 7 alternating maximum deviation points. In Figure~\ref{subfig:error_curve_theta_0.25,0}, there is a shorter peak between -0.5 and 0 which is already in the positive side. This implies that we attain exactly 9 error peaks that come with alternating signs but the magnitude of these error peaks are not uniform. This establishes a good starting point for the Remez method. 

When the correctly identified knot is located in the denominator, Figure~\ref{fig:diff_theta_0.5,0.25} and Figure~\ref{fig:diff_theta_0,0.25} show that the solutions are not optimal as they do not contain the expected number of alternating deviation points in the error curve. Figure~\ref{fig:diff_theta_0.5,0} and Figure~\ref{fig:diff_theta_0,0.5} depict the instance where none of knots are positioned at a correctly identified location. However, Figure~\ref{subfig:error_curve_theta_0,0.5} illustrates that it is more closer to getting an optimal solution than in the case of Figure~\ref{subfig:error_curve_theta_0.5,0} even though the maximum deviation for both cases is identical.

Figure~\ref{fig:diff_theta_(-1/3,1/3)} and Figure~\ref{fig:diff_theta_(1/3,-1/3)} correspond to equidistant knots which divide the domain $[-1,1]$ into 3 equal sized intervals. If the exact locations of the knots are unknown, then it is reasonable to take equidistant knots and ensure that they do not coincide. We consider two scenariors; in the first case, knot in the numerator is located at $-1/3$ and the knot in the denominator is located at $1/3$; in the second case we simply interchange the knot values. Figure~\ref{subfig:error_curve_theta_(-1/3,1/3)} has 7 alternating error peaks whereas Figure~\ref{subfig:error_curve_theta_(-1/3,1/3)} has 9 alternating error peaks although the maximum deviations are exactly the same for both cases.

\begin{figure}
	\begin{center}
		\begin{subfigure}{.49\textwidth}
			\includegraphics[width=\textwidth]{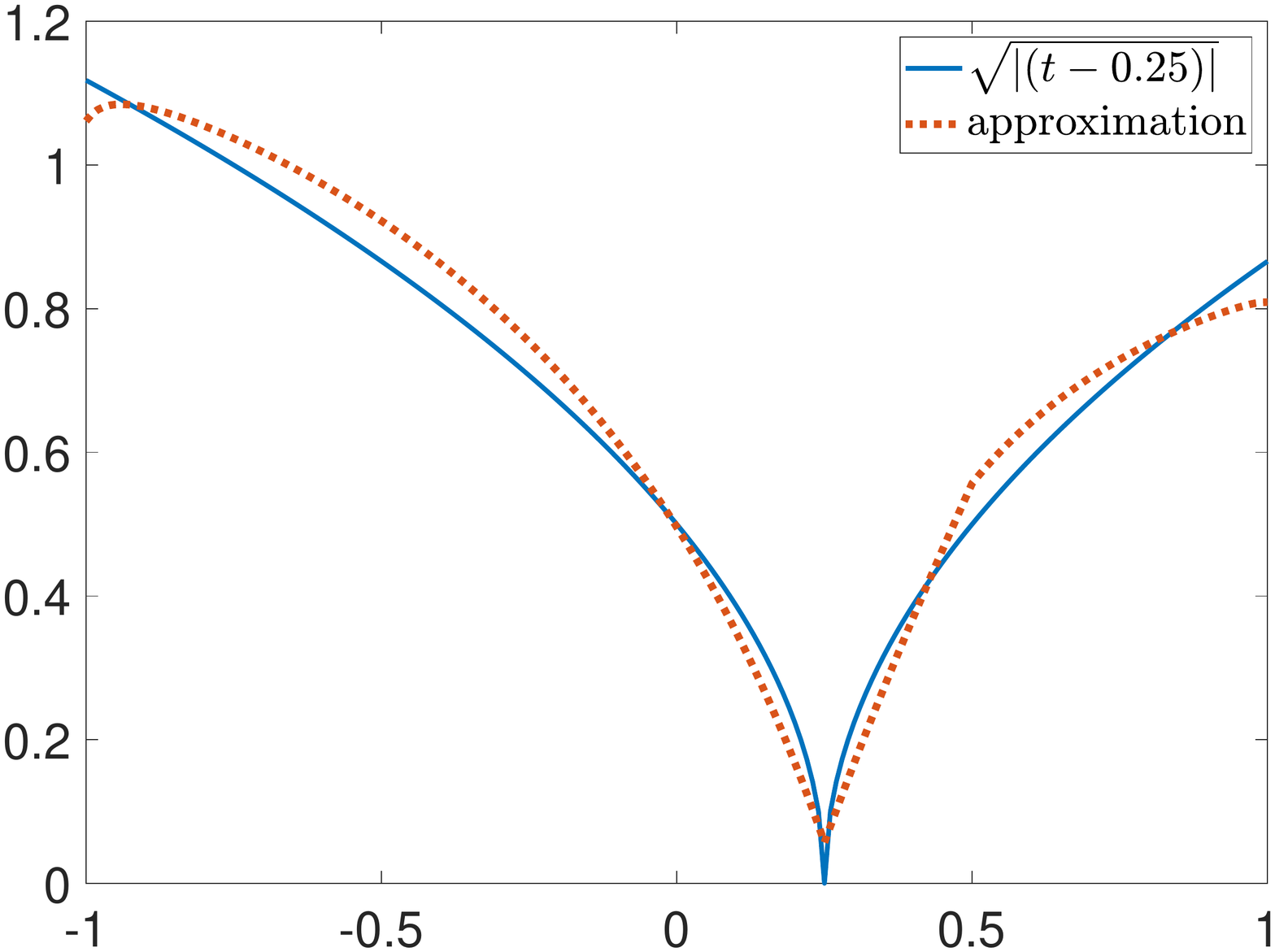}  
			\caption{The function and its approximation}
			\label{subfig:fun_app_theta_0.25,0.5}
		\end{subfigure}
		\begin{subfigure}{.49\textwidth}
			\includegraphics[width=\textwidth]{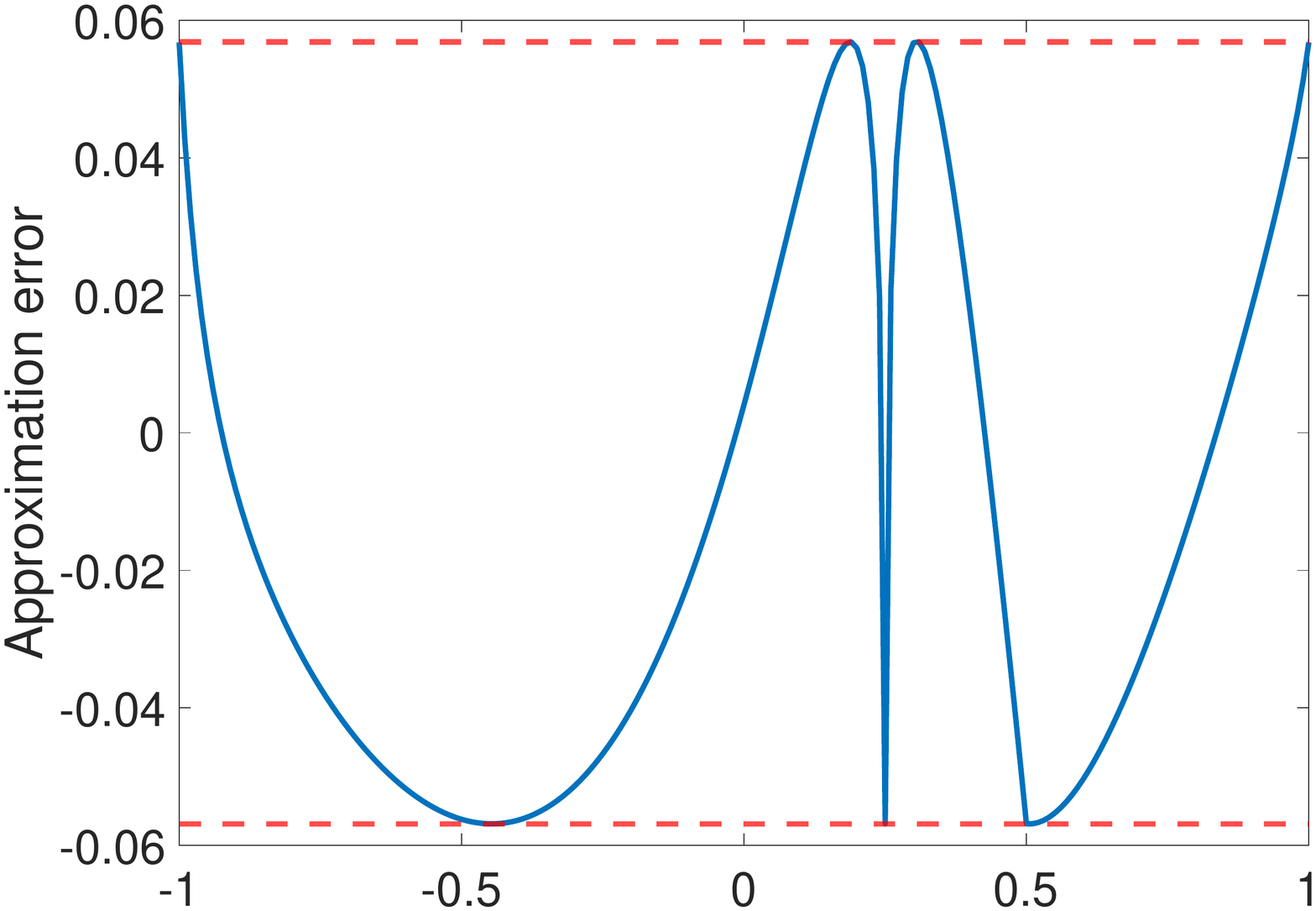}  
			\caption{Error curves}
			\label{subfig:error_curve_theta_0.25,0.5}
		\end{subfigure}
		\caption{The knot in the numerator $\theta_1=0.25$ and the knot in the denominator $\theta_2=0.5$}
		\label{fig:diff_theta_0.25,0.5}
	\end{center}
\end{figure}

\begin{figure}
	\begin{center}
		\begin{subfigure}{.49\textwidth}
			\includegraphics[width=\textwidth]{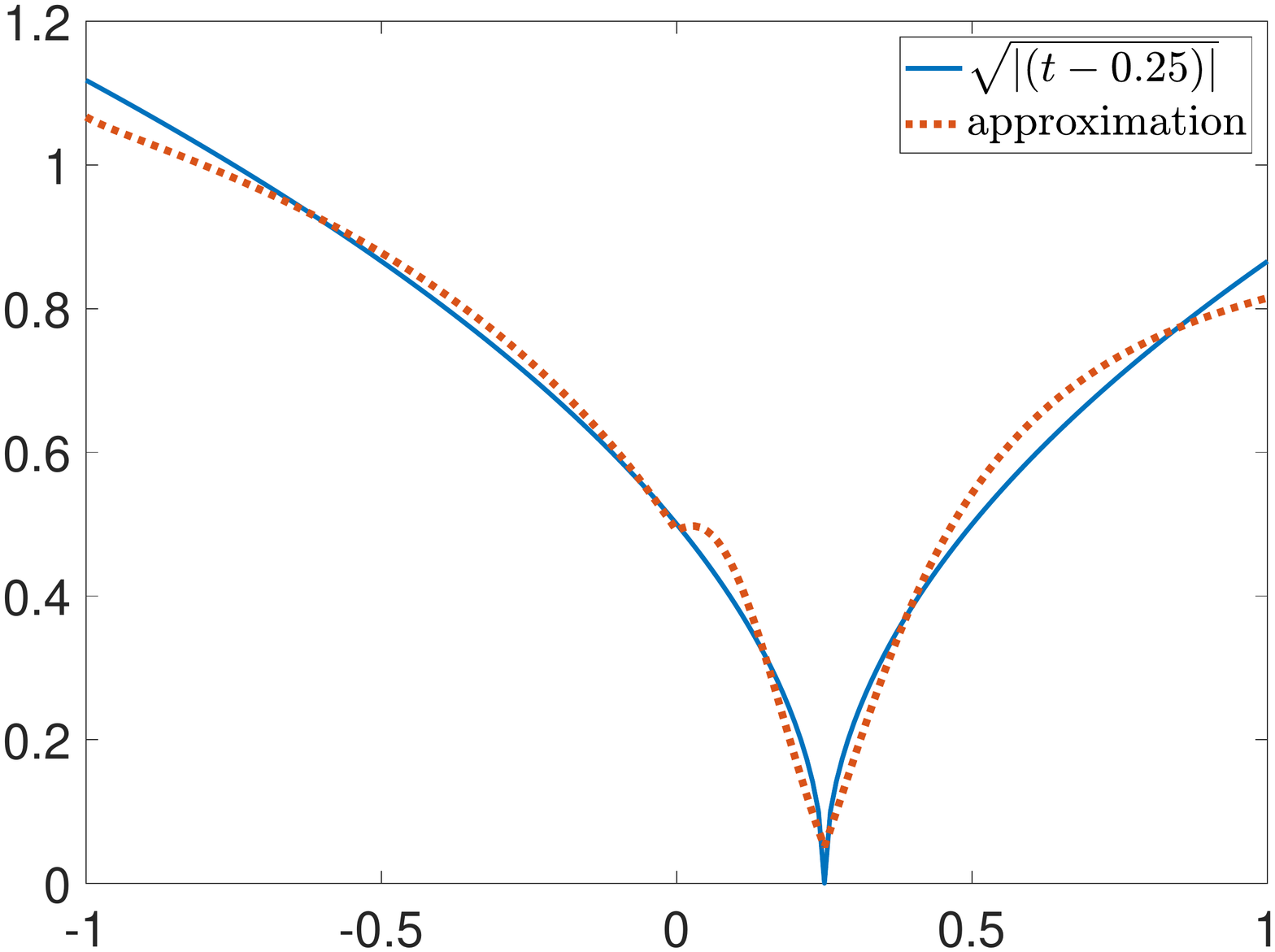}  
			\caption{The function and its approximation}
			\label{subfig:fun_app_theta_0.25,0}
		\end{subfigure}
		\begin{subfigure}{.49\textwidth}
			\includegraphics[width=\textwidth]{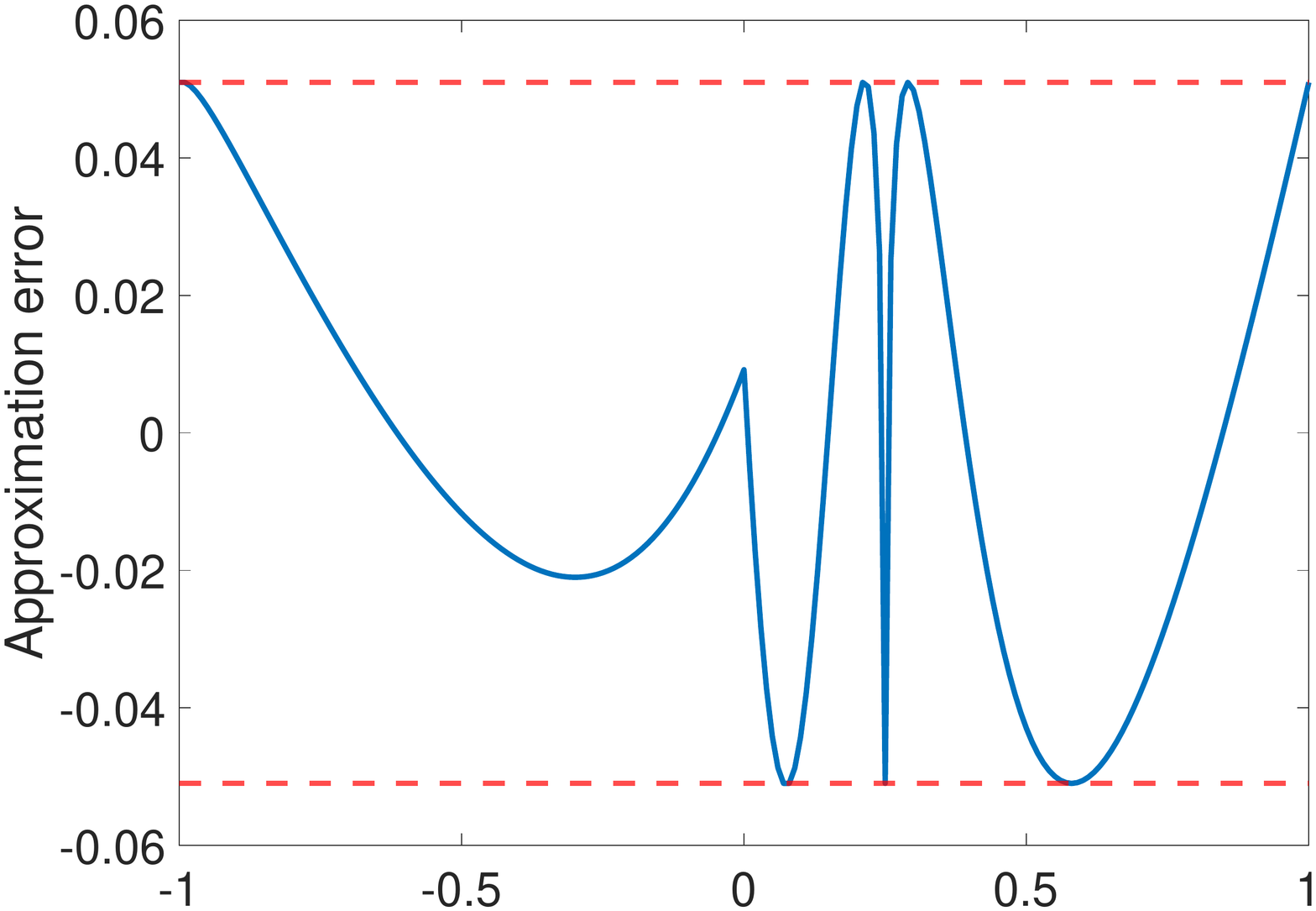}  
			\caption{Error curves}
			\label{subfig:error_curve_theta_0.25,0}
		\end{subfigure}
		\caption{The knot in the numerator $\theta_1=0.25$ and the knot in the denominator $\theta_2=0$}
		\label{fig:diff_theta_0.25,0}
	\end{center}
\end{figure}

\begin{figure}
	\begin{center}
		\begin{subfigure}{.49\textwidth}
			\includegraphics[width=\textwidth]{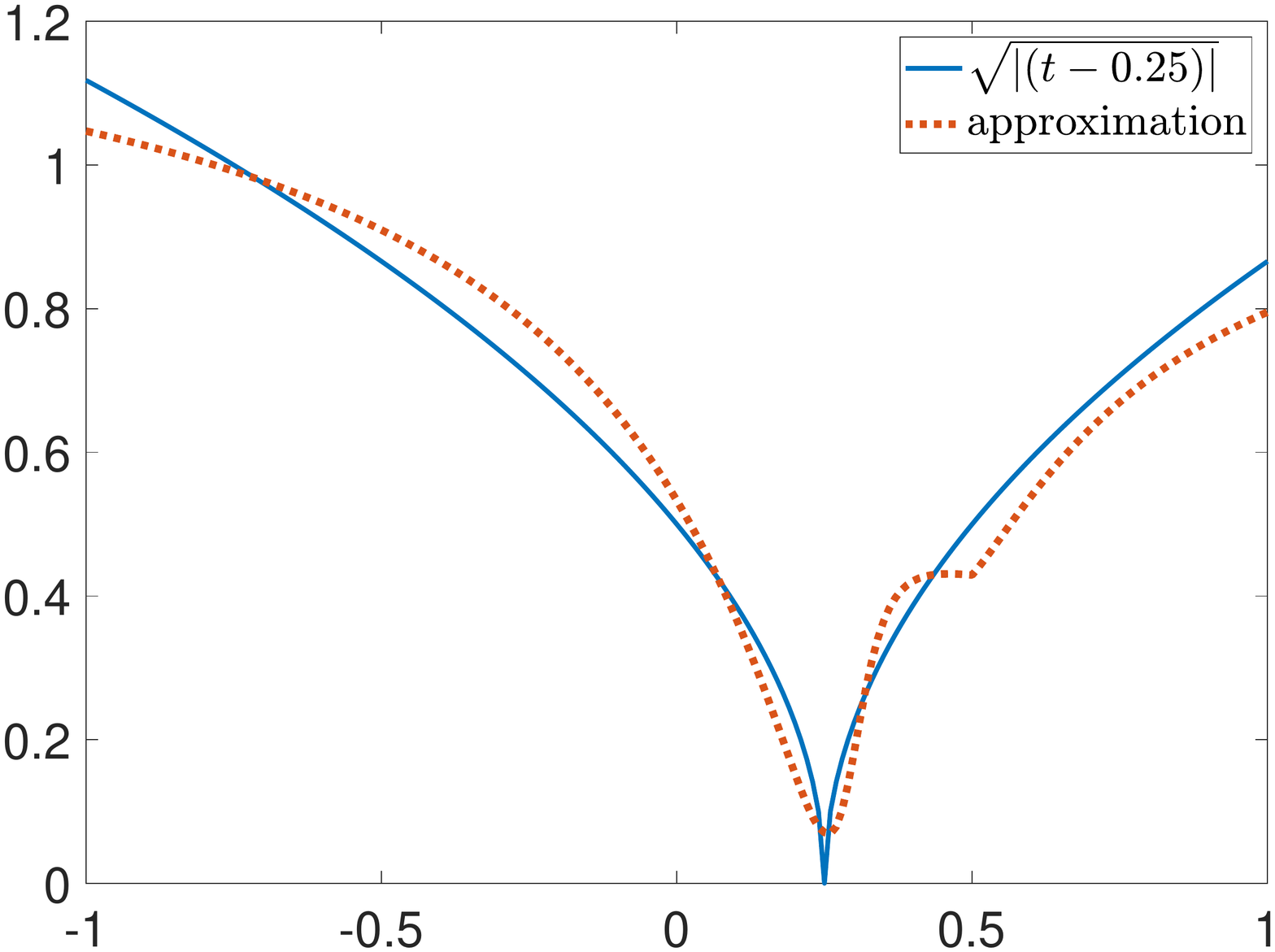}  
			\caption{The function and its approximation}
			\label{subfig:fun_app_theta_0.5,0.25}
		\end{subfigure}
		\begin{subfigure}{.49\textwidth}
			\includegraphics[width=\textwidth]{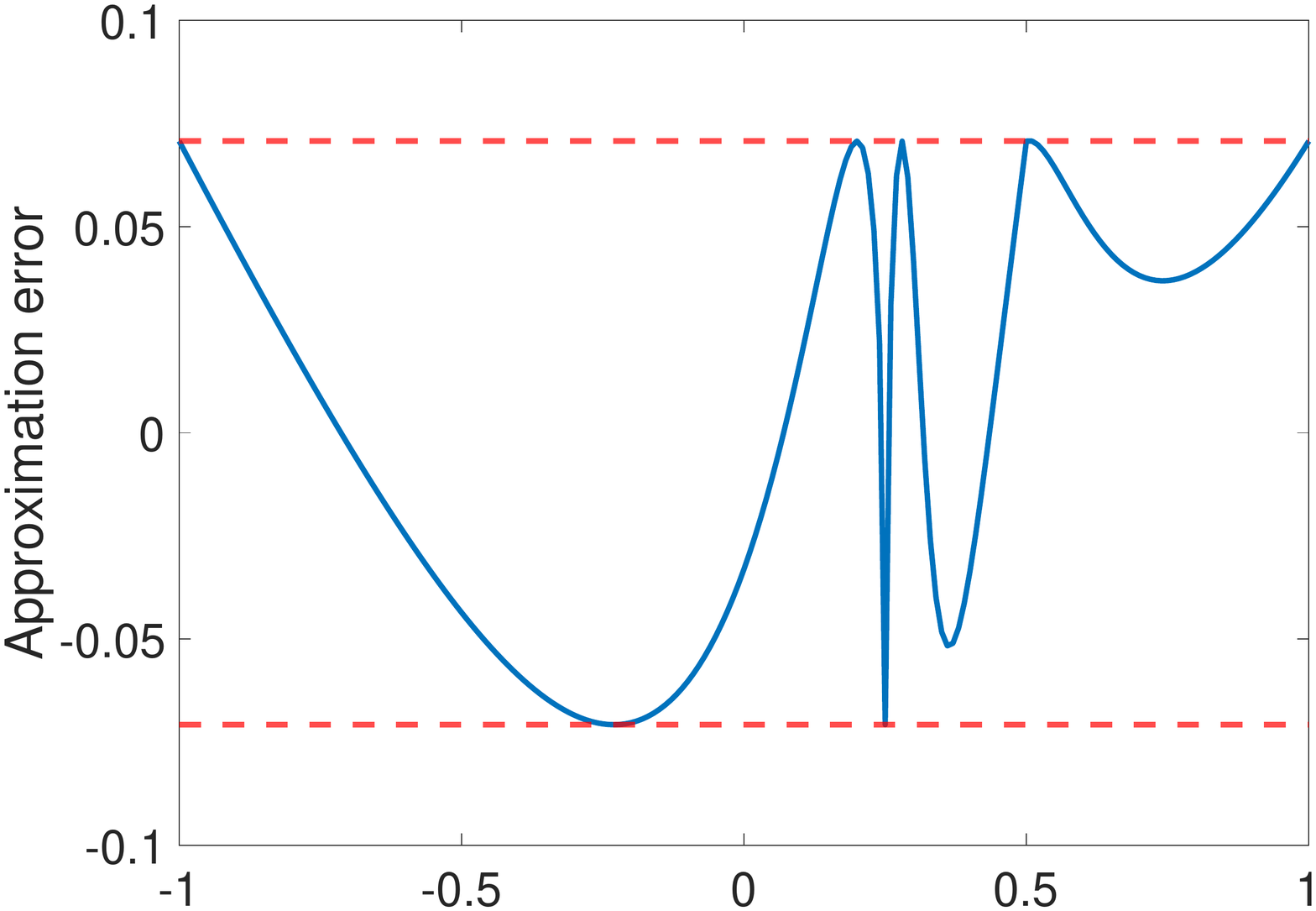}  
			\caption{Error curves}
			\label{subfig:error_curve_theta_0.5,0.25}
		\end{subfigure}
		\caption{The knot in the numerator $\theta_1=0.5$ and the knot in the denominator $\theta_2=0.25$}
		\label{fig:diff_theta_0.5,0.25}
	\end{center}
\end{figure}

\begin{figure}
	\begin{center}
		\begin{subfigure}{.49\textwidth}
			\includegraphics[width=\textwidth]{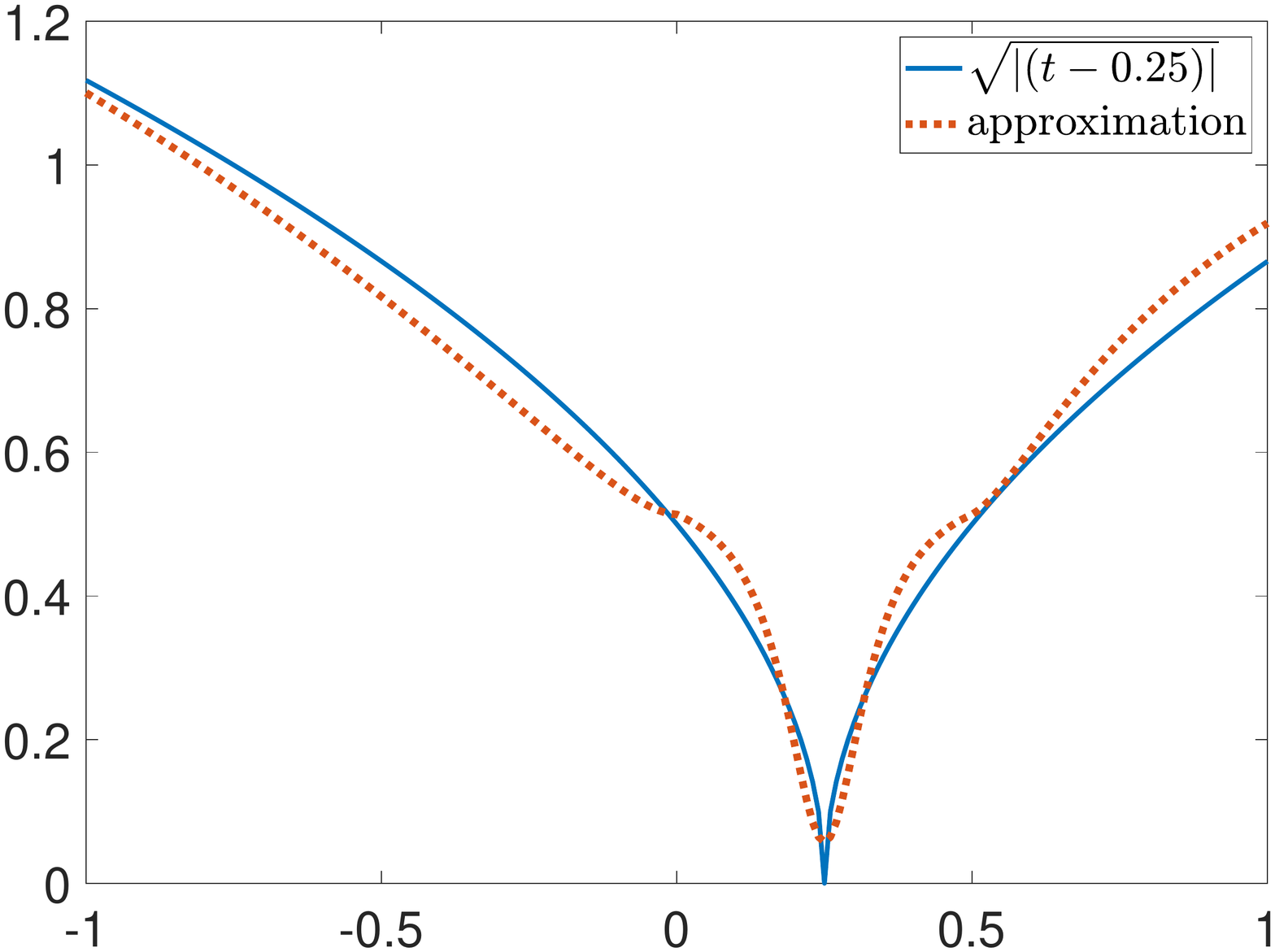}  
			\caption{The function and its approximation}
			\label{subfig:fun_app_theta_0.5,0}
		\end{subfigure}
		\begin{subfigure}{.49\textwidth}
			\includegraphics[width=\textwidth]{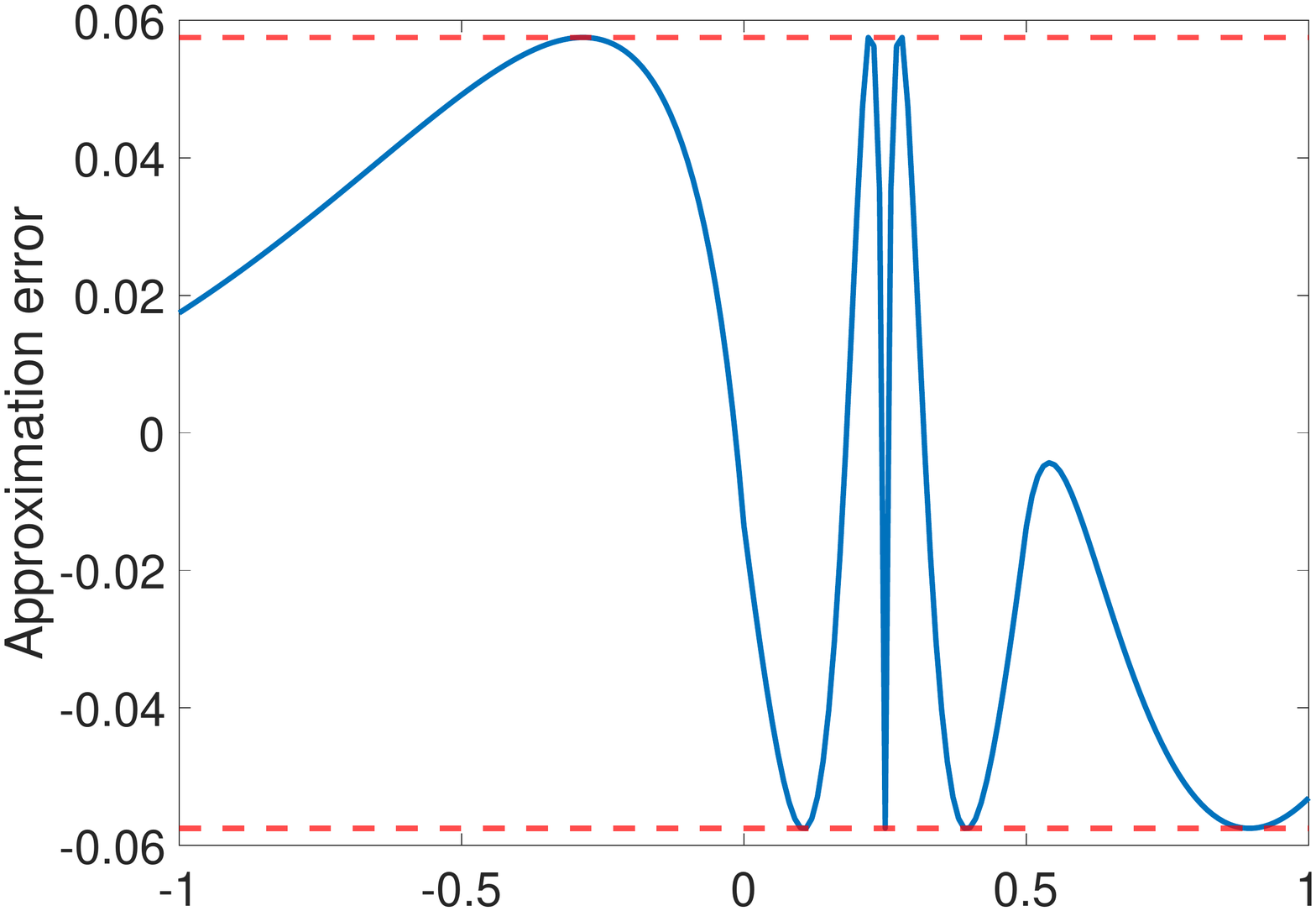}  
			\caption{Error curves}
			\label{subfig:error_curve_theta_0.5,0}
		\end{subfigure}
		\caption{The knot in the numerator $\theta_1=0.5$ and the knot in the denominator $\theta_2=0$}
		\label{fig:diff_theta_0.5,0}
	\end{center}
\end{figure}

\begin{figure}
	\begin{center}
		\begin{subfigure}{.49\textwidth}
			\includegraphics[width=\textwidth]{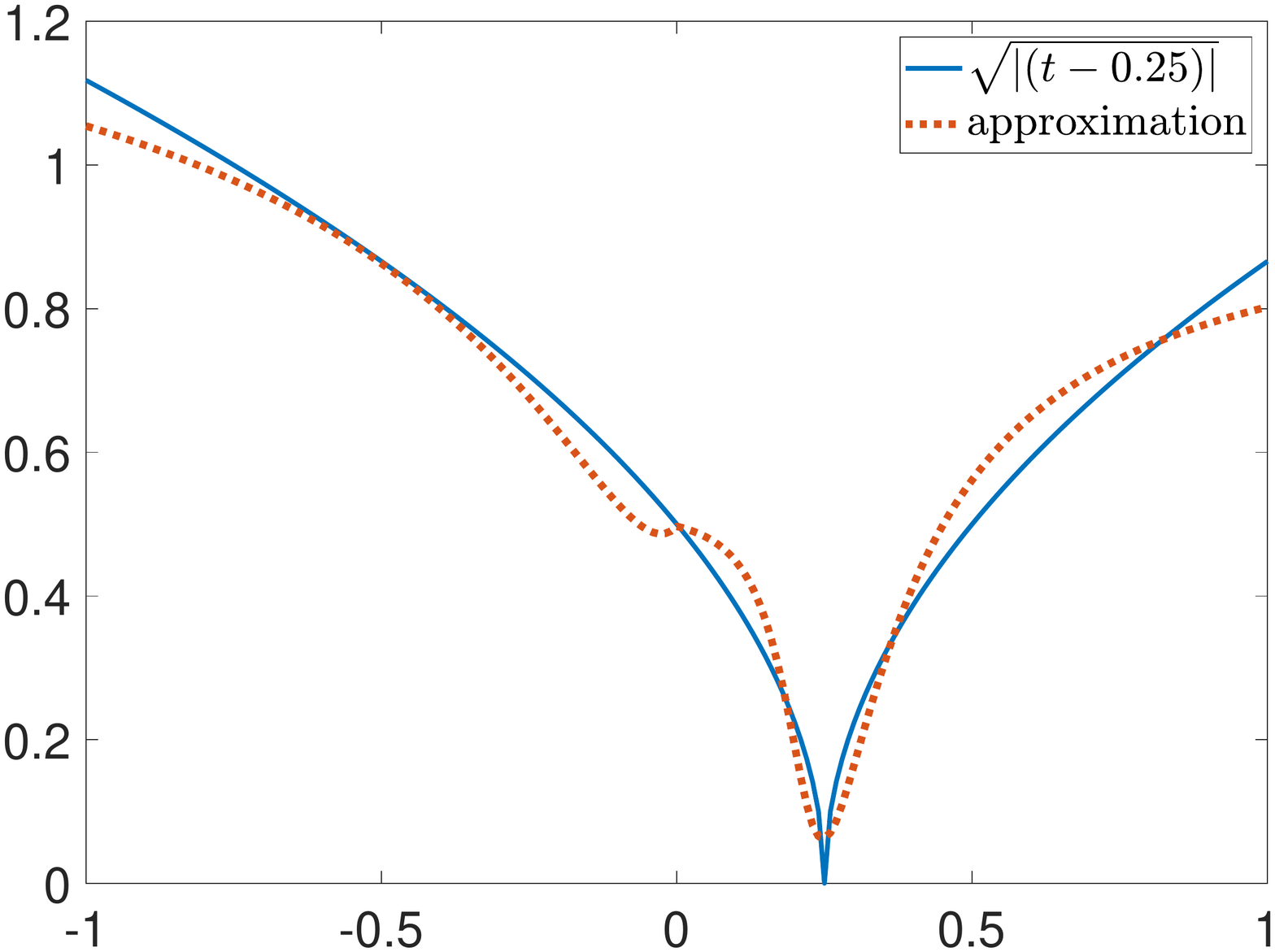}  
			\caption{The function and its approximation}
			\label{subfig:fun_app_theta_0,0.25}
		\end{subfigure}
		\begin{subfigure}{.49\textwidth}
			\includegraphics[width=\textwidth]{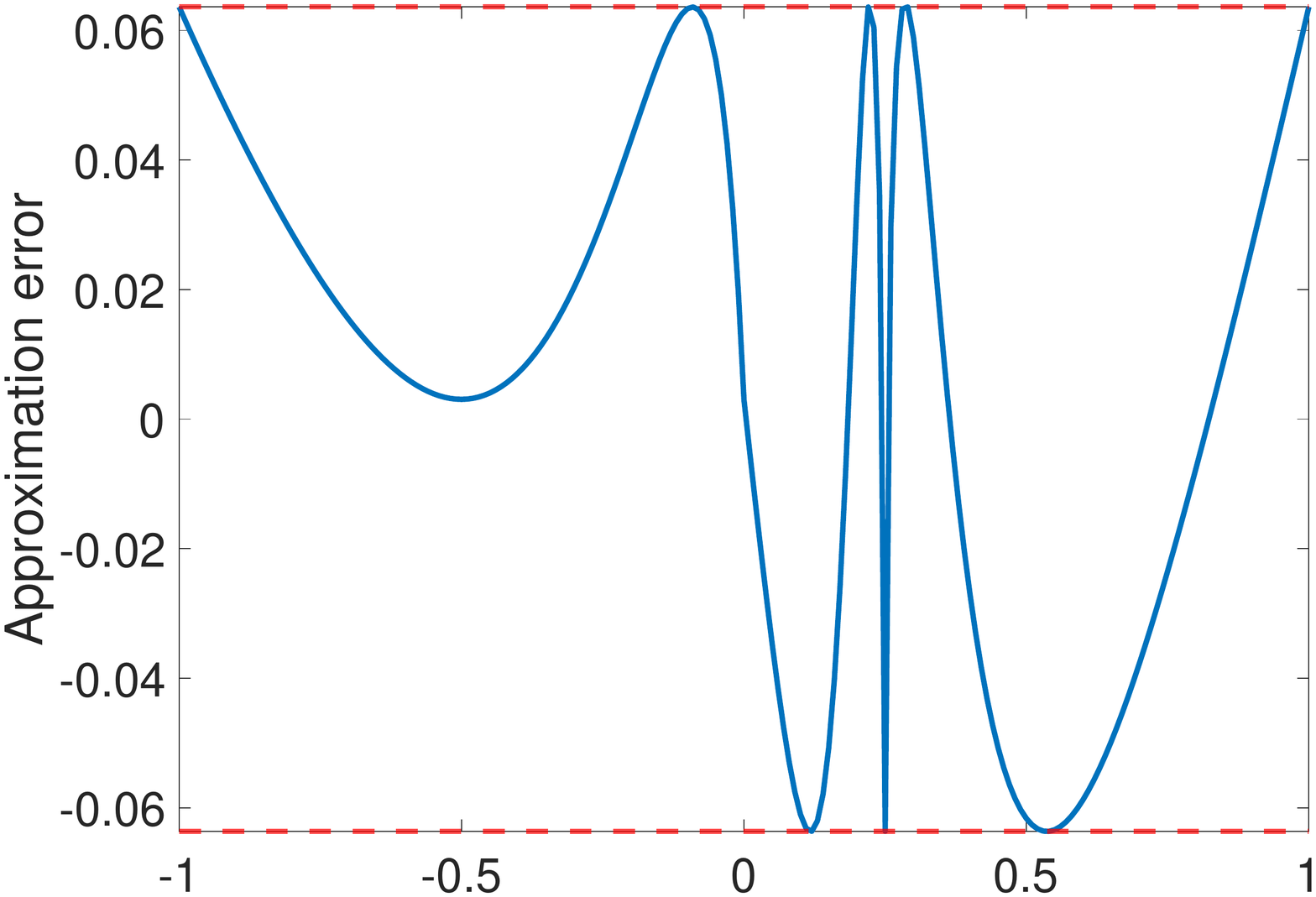}  
			\caption{Error curves}
			\label{subfig:error_curve_theta_0,0.25}
		\end{subfigure}
		\caption{The knot in the numerator $\theta_1=0$ and the knot in the denominator $\theta_2=0.25$}
		\label{fig:diff_theta_0,0.25}
	\end{center}
\end{figure}

\begin{figure}
	\begin{center}
		\begin{subfigure}{.49\textwidth}
			\includegraphics[width=\textwidth]{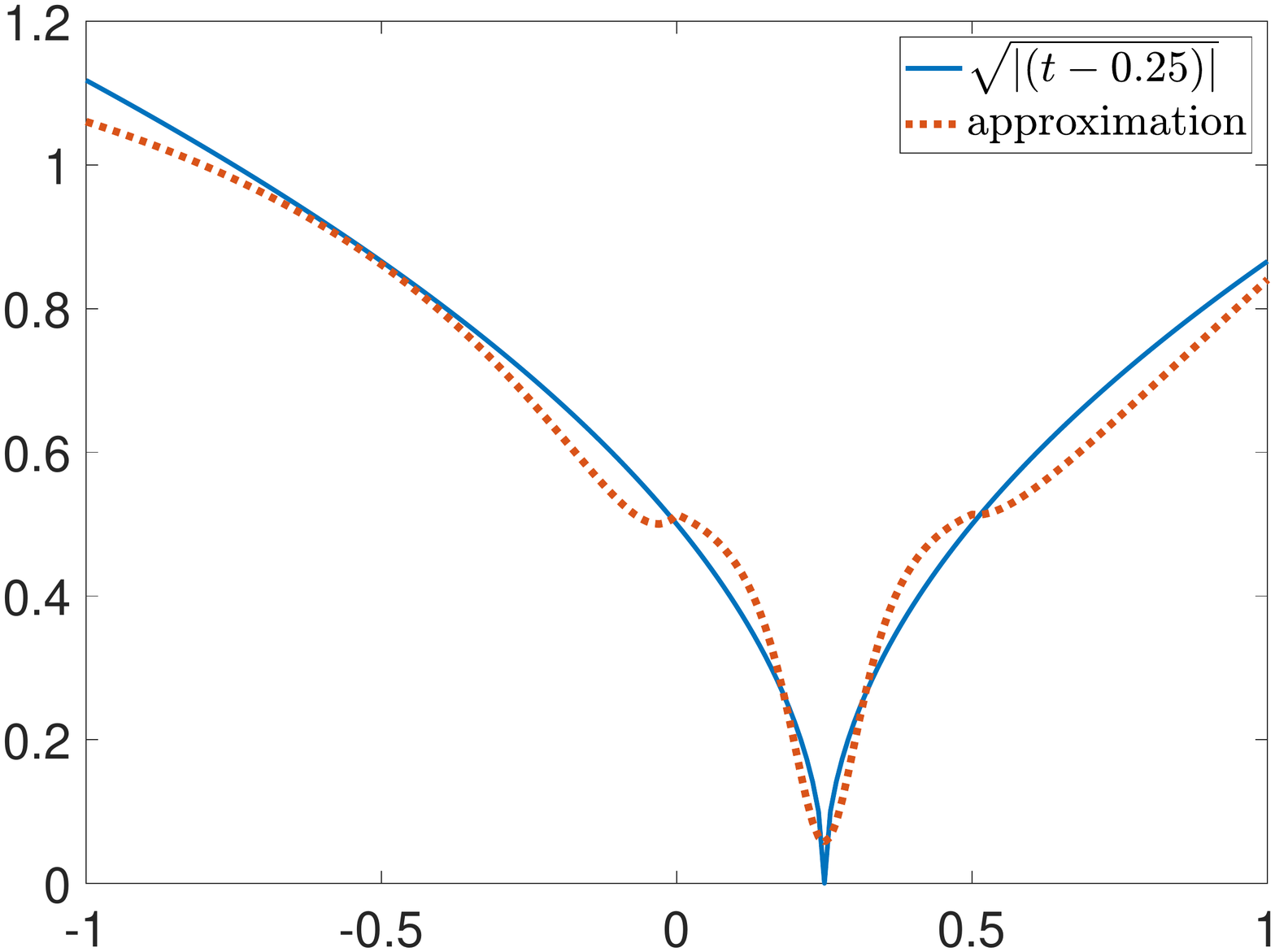}  
			\caption{The function and its approximation}
			\label{subfig:fun_app_theta_0,0.5}
		\end{subfigure}
		\begin{subfigure}{.49\textwidth}
			\includegraphics[width=\textwidth]{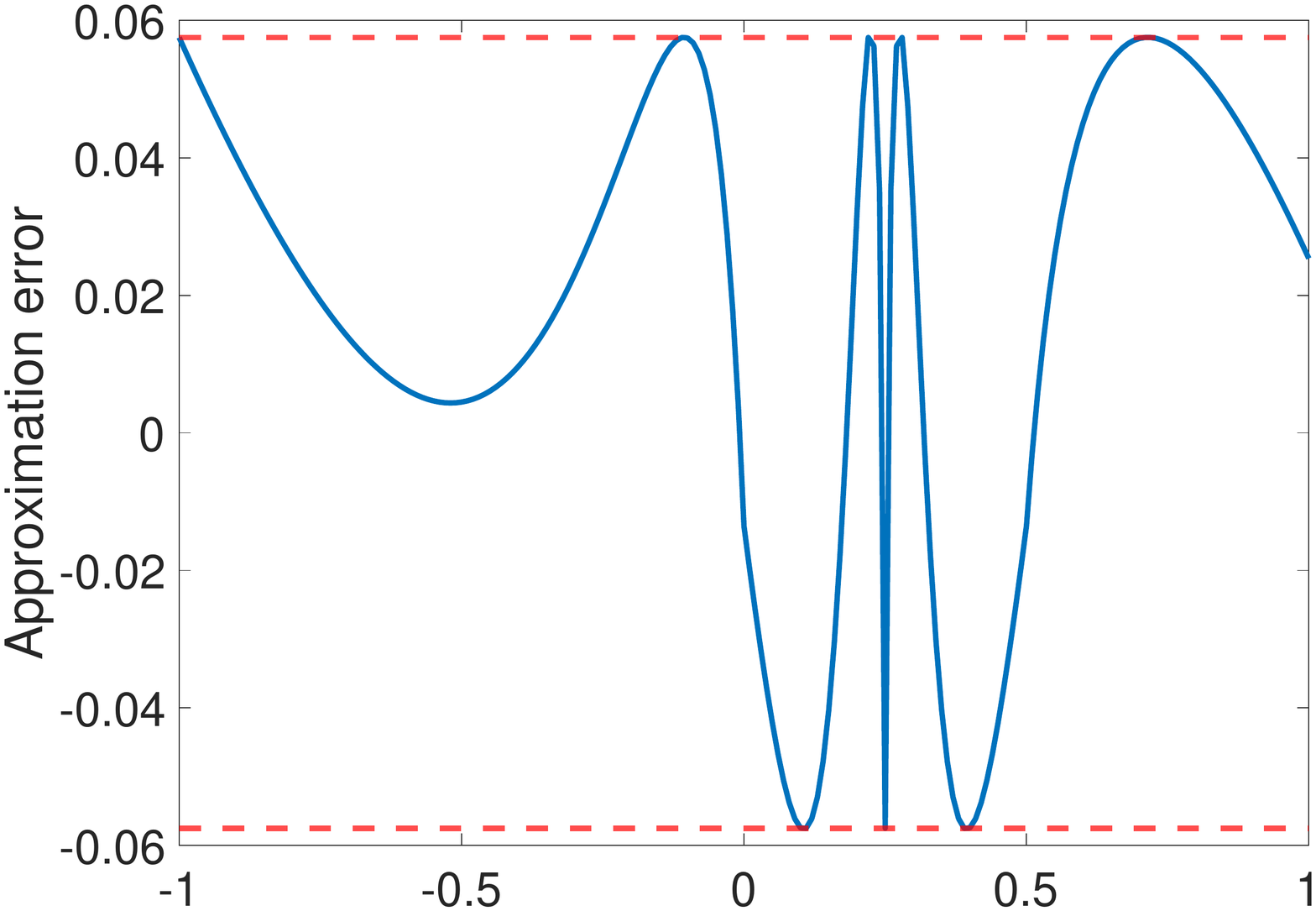}  
			\caption{Error curves}
			\label{subfig:error_curve_theta_0,0.5}
		\end{subfigure}
		\caption{The knot in the numerator $\theta_1=0$ and the knot in the denominator $\theta_2=0.5$}
		\label{fig:diff_theta_0,0.5}
	\end{center}
\end{figure}

\begin{figure}
	\begin{center}
		\begin{subfigure}{.49\textwidth}
			\includegraphics[width=\textwidth]{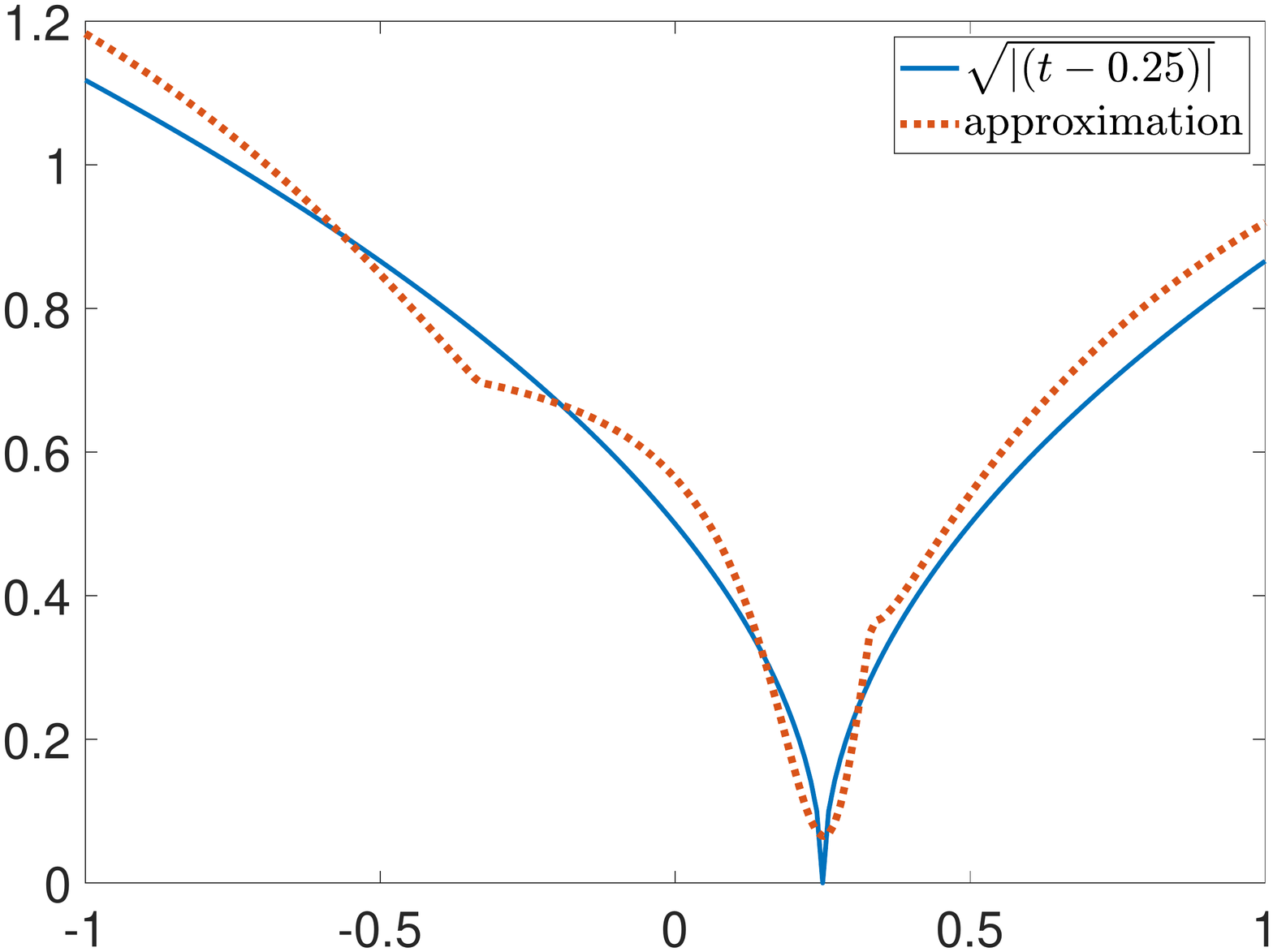}  
			\caption{The function and its approximation}
			\label{subfig:fun_app_theta_(-1/3,1/3)}
		\end{subfigure}
		\begin{subfigure}{.49\textwidth}
			\includegraphics[width=\textwidth]{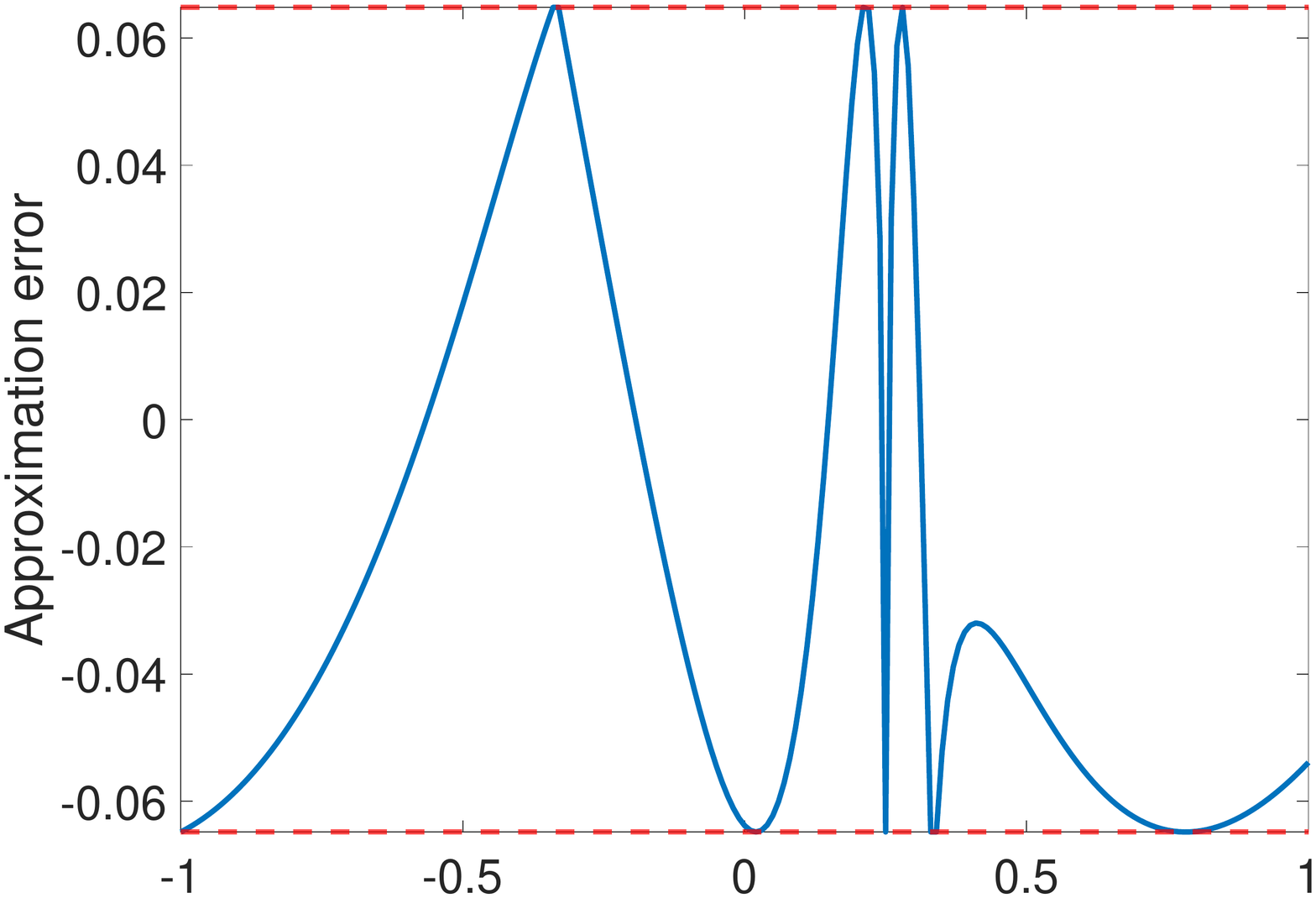}  
			\caption{Error curves}
			\label{subfig:error_curve_theta_(-1/3,1/3)}
		\end{subfigure}
		\caption{The knot in the numerator $\theta_1=-\frac{1}{3}$ and the knot in the denominator $\theta_2=\frac{1}{3}$}
		\label{fig:diff_theta_(-1/3,1/3)}
	\end{center}
\end{figure}

\begin{figure}
	\begin{center}
		\begin{subfigure}{.49\textwidth}
			\includegraphics[width=\textwidth]{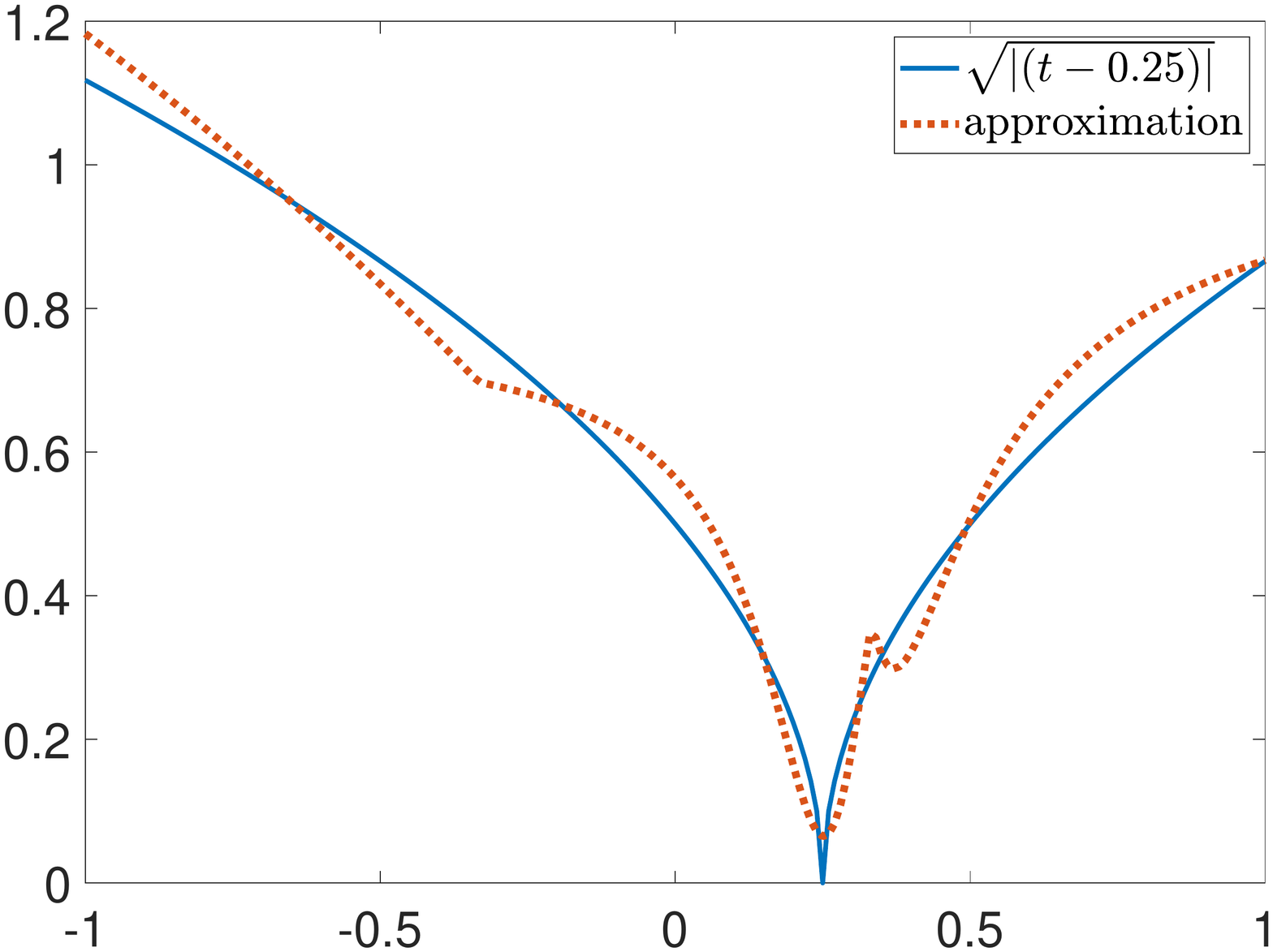}  
			\caption{The function and its approximation}
			\label{subfig:fun_app_theta_(1/3,-1/3)}
		\end{subfigure}
		\begin{subfigure}{.49\textwidth}
			\includegraphics[width=\textwidth]{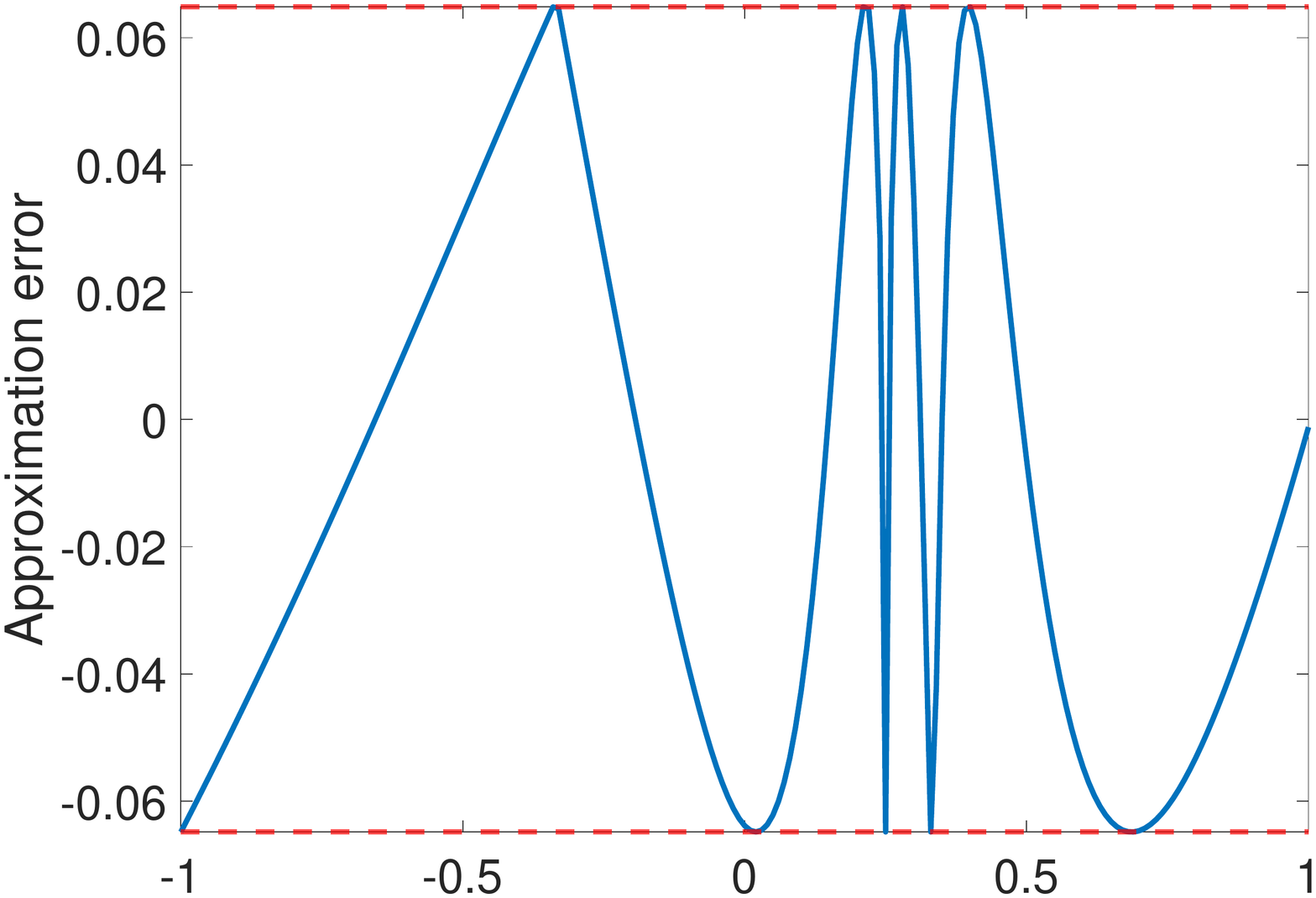}  
			\caption{Error curves}
			\label{subfig:error_curve_theta_(1/3,-1/3)}
		\end{subfigure}
		\caption{The knot in the numerator $\theta_1=\frac{1}{3}$ and the knot in the denominator $\theta_2=-\frac{1}{3}$}
		\label{fig:diff_theta_(1/3,-1/3)}
	\end{center}
\end{figure}

\subsubsection{Second set of experiments}\label{sssec:second_set}

In this set of experiments, the degree of both numerator and denominator of the above approximation is reduced by 1 and the number of knots in each peicewise polynomial is increased by 2 in order to keep the number of parameters unchanged. In particular, the function $f(t)$ is approximated by a ratio of peicewise polynomials of degree 1 with fixed knots. This ratio falls under the class of quasilinear functions leading to Algorithm~\ref{alg:bisection}.

The corresponding optimisation problem is as follows:
\begin{equation*}
	\min_{A,B} \max_{t \in [c,d]} \left| f(t) - S({\bf A},{\bf B},t) \right|,
\end{equation*}
where
${S({\bf A},{\bf B},t)=\frac{a_0 + a_1t + a_2\max\{0,t-\theta_1\} + a_3\max\{0,t-\theta_2\} + a_4\max\{0,t-\theta_3\}}{1 + b_1t + b_2\max\{0,t-\theta_4\} + b_3\max\{0,t-\theta_5\} + b_4\max\{0,t-\theta_6\}}},$\\
$a_i \in A, i=0,\dots,4$ and $b_j \in B, j=1,\dots,4.$

The set of parameters (${\bf A,B}$) consist of $9~(5+4)$ parameters in total to be determined. We consider equidistant knots for the experiments in this section and we start by dividing the domain into 4 equal sized intervals. This helps to reduce the overlapping between the knots in the numerator and the denominator and it provides more flexibility for the knots. Hence the number of parameters to be determined will be the same as before.

Assume $\theta_1$, $\theta_2$ and $\theta_3$ (knots in the numerator) are distinct and they are equidistant. The knots in the denominator, $\theta_4$, $\theta_5$ and $\theta_6$ coincide with the knots in the numerator. 
Therefore, we take, 
$$\theta_1 = -0.5 = \theta_4,$$
$$\theta_2 = 0 = \theta_5,$$ 
$$\theta_3 = 0.5 = \theta_6$$.

\begin{figure}
	\begin{center}
		\begin{subfigure}{.49\textwidth}
			\includegraphics[width=\textwidth]{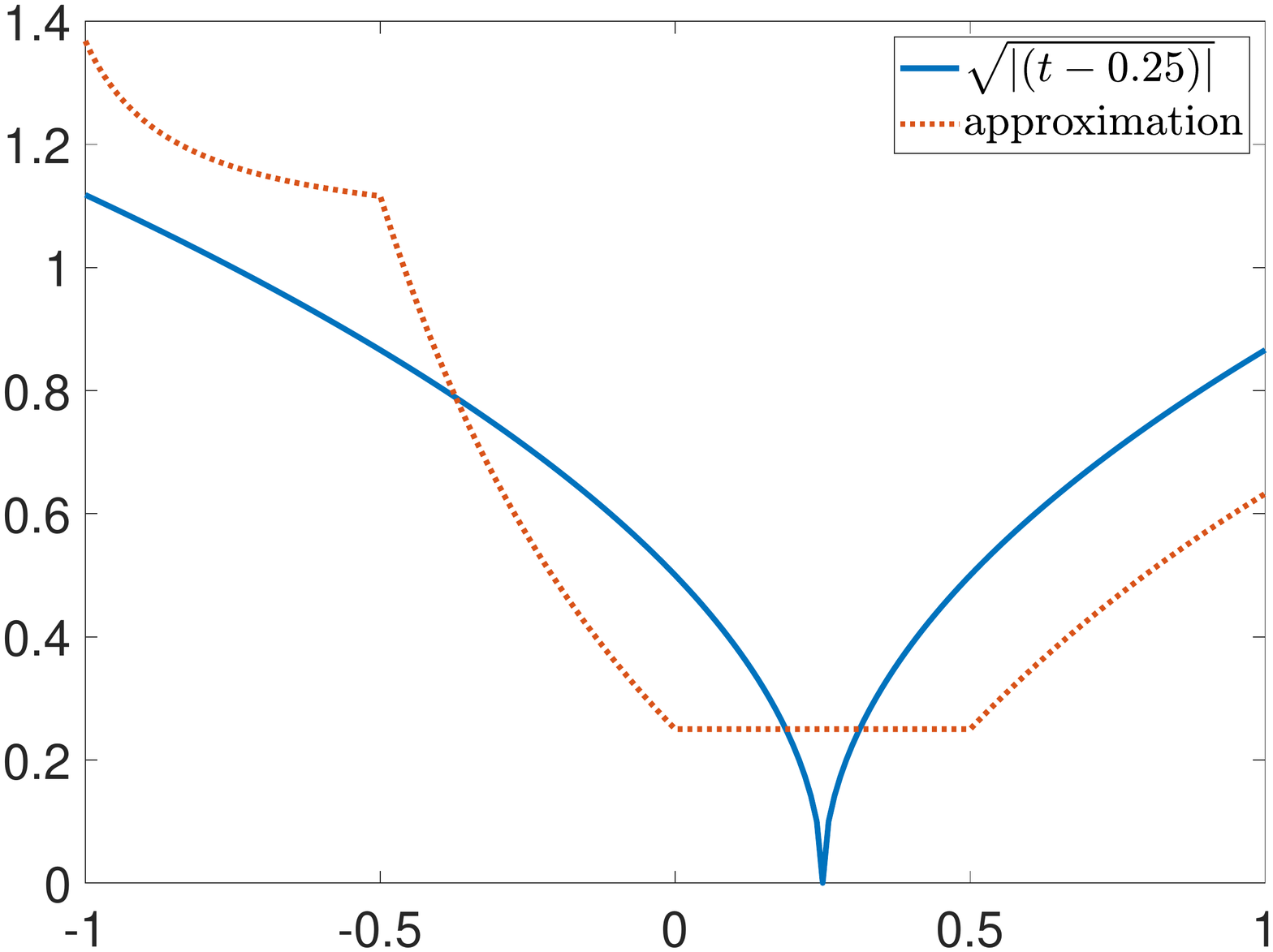}  
			\caption{The function and its approximation}
			\label{subfig:fun_app_4_intervals}
		\end{subfigure}
		\begin{subfigure}{.49\textwidth}
			\includegraphics[width=\textwidth]{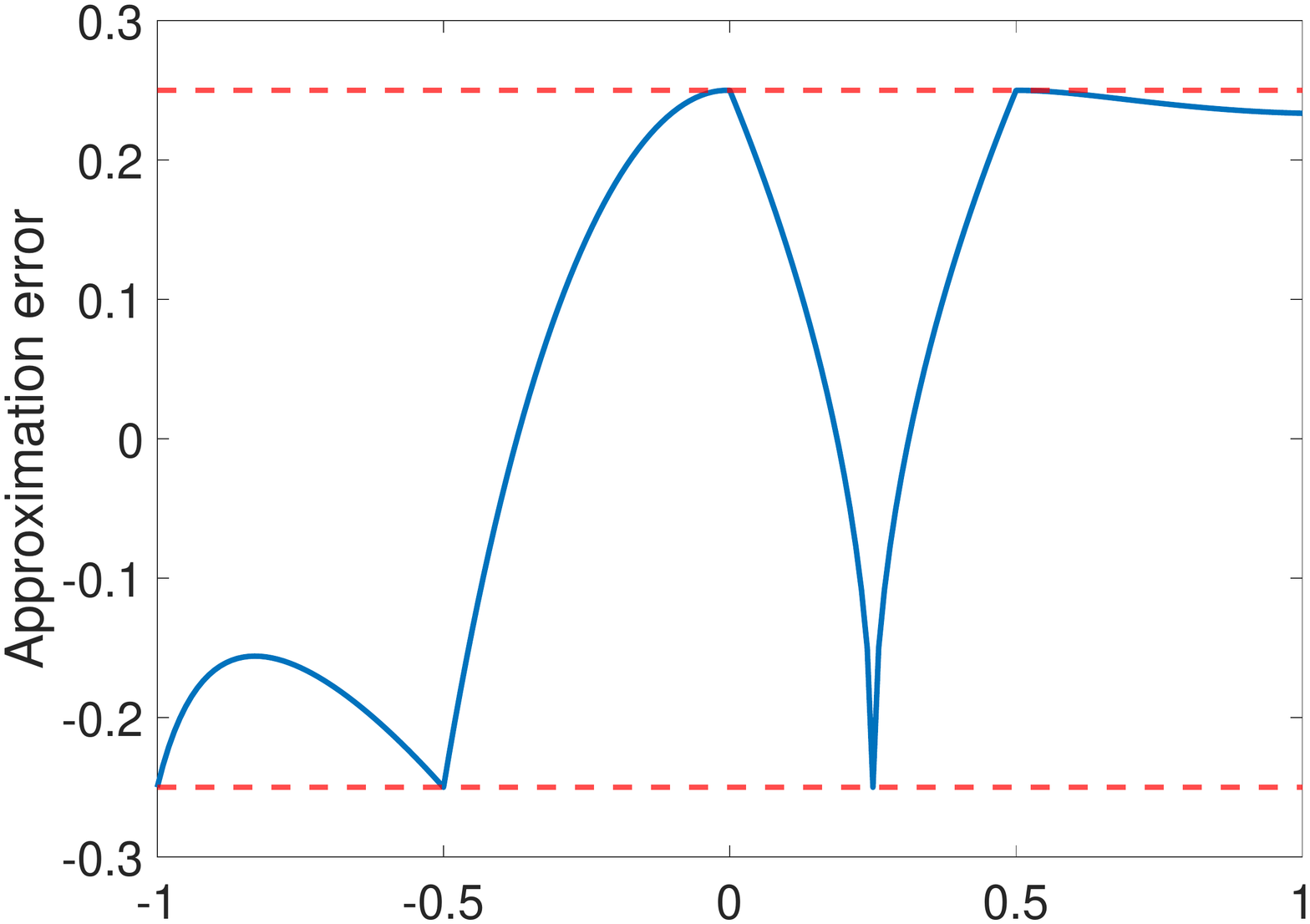}  
			\caption{Error curves}
			\label{subfig:error_curve_4_intervals}
		\end{subfigure}
		\caption{$\theta_1 = -0.5 = \theta_4$, $\theta_2 = 0 = \theta_5$, $\theta_3 = 0.5 = \theta_6.$}
		\label{fig:diff_theta_4_intervals)}
	\end{center}
\end{figure}

Figure~\ref{subfig:fun_app_4_intervals} clearly indicates that the approximation is made up of 4 small rigid pieces of functions and Figure~\ref{subfig:error_curve_4_intervals} only contain 4 alternating error peaks with maximum deviation being around 0.25. One can expect to get an accurate approximation when the number of knots is increased. But if none of the knots are positioned at a correctly identified location or at least somewhere close to an original (known) knot point, it is very unlikely to get an optimal solution.

\begin{figure}
	\begin{center}
		\begin{subfigure}{.49\textwidth}
			\includegraphics[width=\textwidth]{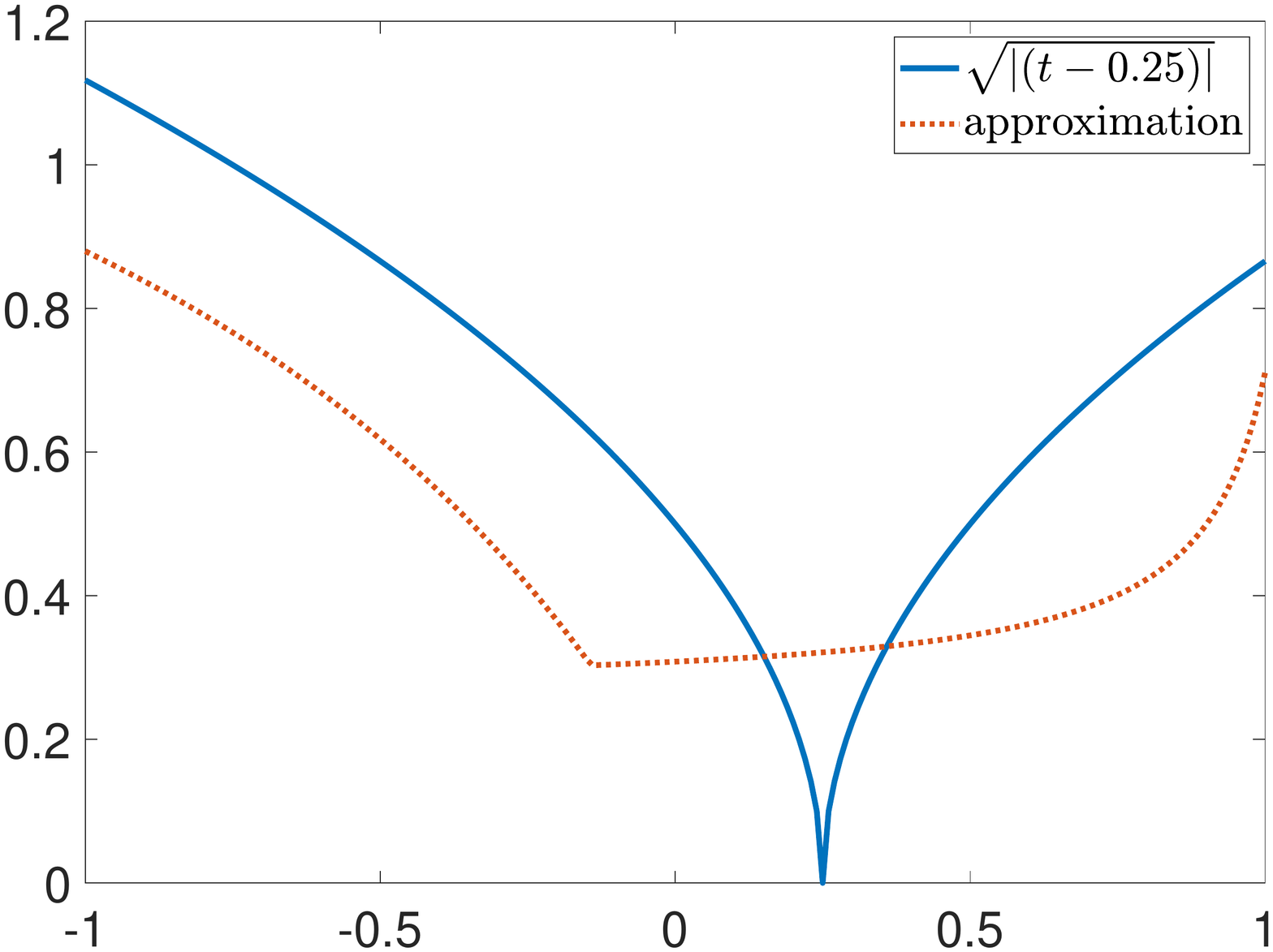}  
			\caption{The function and its approximation}
			\label{subfig:fun_app_7_intervals}
		\end{subfigure}
		\begin{subfigure}{.49\textwidth}
			\includegraphics[width=\textwidth]{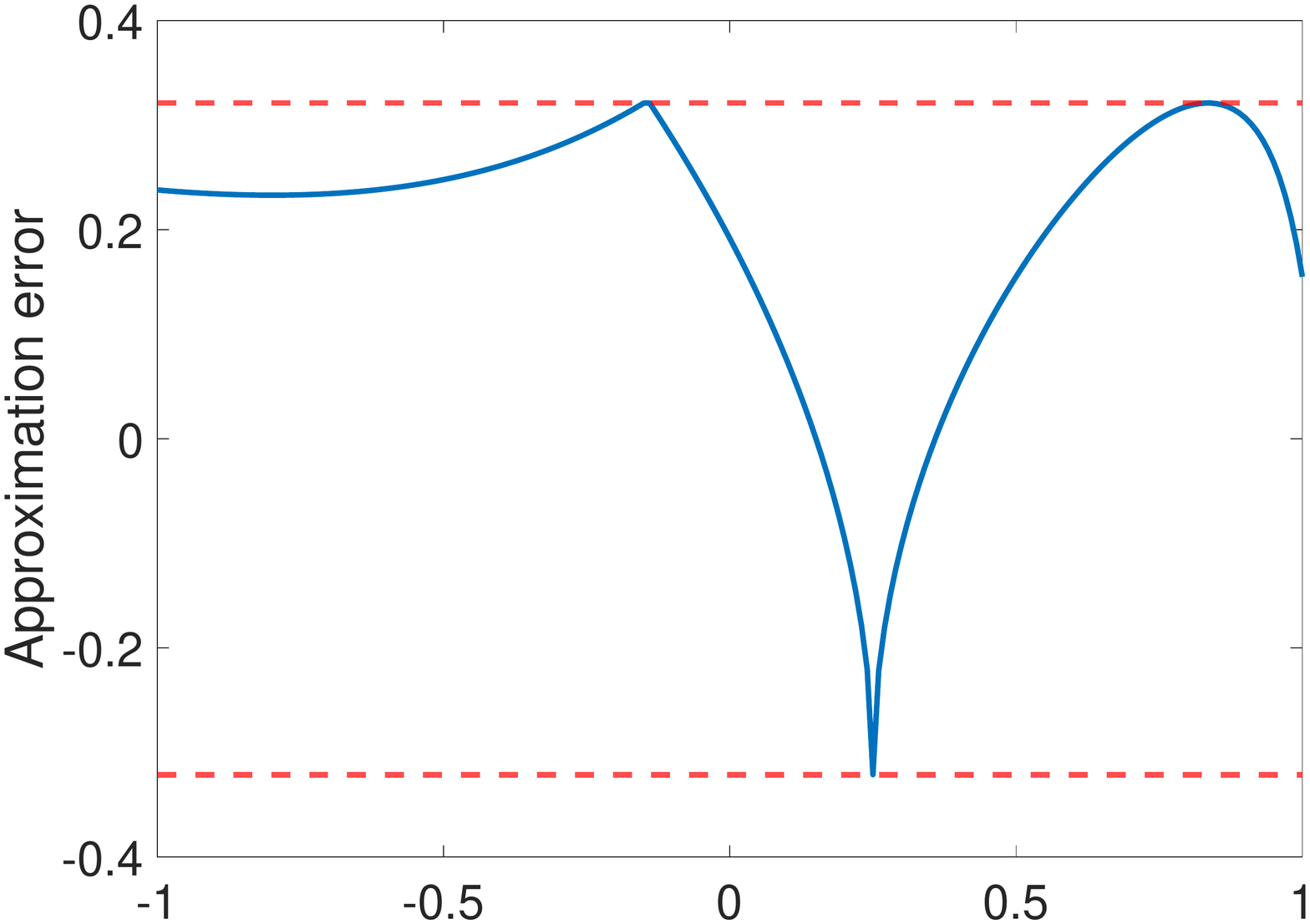}  
			\caption{Error curves}
			\label{subfig:error_curve_7_intervals}
		\end{subfigure}
		\caption{$\theta_1=-5/7, \theta_2=-3/7, \theta_3=-1/7, \theta_4=1/7, \theta_5=3/7, \theta_6=5/7.$}
		\label{fig:diff_theta_7_intervals)}
	\end{center}
\end{figure}

Another interesting aspect is shown in Figure~\ref{fig:diff_theta_7_intervals)} where all 6 knots in the approximation are distributed equidistantly by dividing the domain into 7 equal sized intervals. 

\begin{equation*}
	\begin{split}
		\theta_1 = -\frac{5}{7},\\
		\theta_4 = \frac{1}{7}, 
	\end{split}
	\quad
	\begin{split}
		\theta_2 = -\frac{3}{7},\\
		\theta_5 = \frac{3}{7},
	\end{split}
	\quad
	\begin{split}
		\theta_3 = -\frac{1}{7},\\
		\theta_6 = \frac{5}{7}.
	\end{split}
\end{equation*}

There is no overlapping between the knots in the numerator and the denominator. However, the error curve contains only 3 alternating peaks and the maximum deviation is much higher than any other cases that we considered earlier.

\section{Conclusions and future work}\label{sec:conclusions}

In this paper we demonstrated that the linear inequality method developed decades ago for rational and generalised rational Chebyshev approximation is the direct application of the bisection method for quasiconvex optimisation. This correspondence is not surprising, since it is a well-known fact that the optimisation problems in generalised rational chebyshev approximation is quasiconvex. However, this observation leads to an extension of this method to a broader class of functions.

In the future, we are planning to continue our work in the direction of identifying the correspondence between generalised rational Chebyshev approximation methods and general optimisation methods for quasiconvex and psedoconvex problems. Another interesting research direction is to investigate the extension of generalised Chebyshev approximation to multivariate settings, where the application of general quasiconvex optimisation methods look very promising.     

\section*{Acknowledgement}
This research was supported by the Australian Research Council (ARC),  Solving hard Chebyshev approximation problems through nonsmooth analysis (Discovery Project DP180100602).


\end{document}